\documentclass[11pt]{article}
\usepackage[margin=1in]{geometry}
\usepackage{amsmath,amssymb,amsthm}
\usepackage{graphicx}
\usepackage{booktabs}
\usepackage{enumitem}
\usepackage{xcolor}
\usepackage[colorlinks=true,linkcolor=blue!50!black,citecolor=blue!50!black,urlcolor=blue!50!black]{hyperref}
\graphicspath{{./}{./figs/}}

\usepackage{algorithm}
\usepackage{algpseudocode}
\usepackage{fancyhdr}
\usepackage{placeins}

\newcommand{\Hcal}{\mathcal H}

\newcommand{\R}{\mathbb{R}}
\newcommand{\Sph}{\mathbb{S}}
\newcommand{\tr}{\operatorname{tr}}
\newcommand{\dd}{\,d\mu}
\newcommand{\Span}{\operatorname{span}}

\newtheorem{theorem}{Theorem}
\newtheorem{proposition}[theorem]{Proposition}
\newtheorem{lemma}[theorem]{Lemma}
\newtheorem{corollary}[theorem]{Corollary}

\theoremstyle{definition}

\theoremstyle{remark}
\newtheorem{remark}[theorem]{Remark}

\fancypagestyle{coverpage}{
    \fancyhf{} 
    
    \fancyfoot[L]{%
        \parbox[b]{0.85\textwidth}{%
            \footnotesize \textit{Generative AI was used for copy editing and language polishing.}%
        }%
    }
    
    \fancyfoot[R]{\thepage}
}

\title{Spectral Landmarks: \\ Alignment of Landmark Kernel Translates \\ with the Leading Eigenfunctions of a Kernel Operator}
\author{Jeanie Schreiber}
\date{}

\begin{document}
\let\oldfootnoterule\footnoterule
\renewcommand{\footnoterule}{}
\maketitle
\thispagestyle{coverpage}

\begin{abstract}
We study the problem of choosing landmark points $Z=\{z_1,\dots,z_m\}$ so that the span of the kernel translates $k(\cdot,z_j)$ is as close as possible to the leading $s$-dimensional eigenspace of a positive-definite kernel. The natural measure is the Hilbert-Schmidt operator discrepancy $\tfrac12\|\Pi_Z-Q_s\|_{HS}^2$ between the reproducing-kernel-Hilbert-space (RKHS) projectors onto the two subspaces. We show that $F(Z)=\tfrac12\|\Pi_Z-Q_s\|_{HS}^2 -\frac{m-s}{2}$ has a closed form in terms of the truncated Mercer kernel and relate it to the classical integrated power function (Nystr\"om trace error) through a two-sided sandwich bound that holds for every $m\ge s$. Minimizing $F$ by gradient descent on the landmark coordinates produces landmark sets that are \emph{evenly spaced}. Three theoretical results support this observation. (i) For $m=s$, writing $F$ through the tangents of the principal angles yields an upper bound minimized precisely by \emph{tight sampling sets}, i.e.\ configurations whose eigenfunction sampling matrix is a scaled isometry. On homogeneous spaces these are exactly the minimal equal-weight cubature rules for products of leading eigenfunctions, and their natural extension to every $m\ge s$ consists of the configurations whose sampled eigenfunctions form a unit-norm tight frame. (ii) On homogeneous spaces, a symmetry principle shows that every orbit of a finite isometry group whose stabilizers fix no tangent vector is a critical point. (iii) On the circle, a closed-form aliasing formula gives the value of $F$ at equispaced landmarks for every $m\ge s$. Experiments on $\Sph^1$ and $\Sph^2$ agree with the theory. On the circle, every random start converges to equispaced landmarks, and on the sphere the optimizer recovers the regular tetrahedron, octahedron, and icosahedron from random starts. Finally, we analyze the numerical behavior of the problem. \end{abstract}


\tableofcontents

\section{Introduction}
Low-rank kernel approximations aim to circumvent the $\mathcal{O}(n^2)$ memory and $\mathcal{O}(n^3)$ factoring cost of large $n\times n$ Gram matrices by focusing on strategically chosen subsets of kernel entries. For symmetric positive definite kernels, this challenge manifests in identifying a subset that preserves spectral structure while minimizing approximation error.

The Nystr\"om method exemplifies this approach by selecting $m\ll n$ \emph{landmarks} $Z=\{z_1,\dots,z_m\}$ and
replacing the kernel $K$ by its projection $\widehat K$ onto the span of the translates
$k(\cdot,z_j)$
\cite{WilliamsSeeger2001,DrineasMahoney2005,GittensMahoney2016}. The goal is a
low-rank surrogate of $K$ that is cheap to store and apply, yet close to
$K$ in a norm that matters for the downstream task. The quality of the approximation hinges critically on the selection of the landmark points, a combinatorial optimization problem with no universally efficient solution.

A natural way forward is to measure the residual $K-\widehat K$ directly.
Its diagonal is the squared \emph{power function} $P_Z^2$ of
scattered-data approximation, which is the squared RKHS distance from
$k(\cdot,x)$ to the span of the translates. Integrating $P_Z^2$ gives the trace error
of the Nystr\"om approximation, so minimizing $\int P_Z^2$, or $\sup P_Z$,
or running greedy selection are all residual-based criteria for
landmarks which control the point-wise
error of kernel interpolation
\cite{Schaback1995,Wendland2005}.

Many uses of a kernel depend on its leading $s$-dimensional eigenspace $V_s$ of the integral operator, which carries the principal components of the feature map and, after the appropriate degree normalization of the kernel, the leading diffusion-map coordinates
\cite{CoifmanLafon2006}. However, a landmark set can have a small trace error and
still span a subspace that is rotated away from $V_s$, in particular when
the spectrum has no clear gap at $s$. This paper therefore studies a different, explicitly spectral criterion.
Rather than asking that the translates approximate $K$ well \emph{on
average}, we ask that their span align with $V_s$ explicitly. The resulting
\emph{spectral alignment} objective is a Grassmannian distance between two
subspaces of the RKHS. It is smooth in the landmark positions and inexpensive to
evaluate. After a one-time computation of the top $s$ eigenpairs, each evaluation
needs the $m\times N$ cross-kernel matrix between landmarks and data and a small $m\times m$ Cholesky factorization. Its minimizers turn out to be \emph{evenly spaced},
which we support with theory (Section~\ref{sec:evenness}) and with numerical
experiments on the circle and the sphere
(Sections~\ref{sec:circle}--\ref{sec:sphere}).

We begin by deriving a closed form for the spectral subspace alignment in terms of the truncated and tail Mercer kernels, and we show that it equals the sum of squared sines of the principal angles between the span of the kernel translates and the leading eigenspace $V_s$ (Section~\ref{sec:discrepancy}). This form has a useful dual reading as measuring both the mass of $V_s$ that leaks out of the landmark space and the mass of the kernel tail that leaks into it. We then connect the alignment functional to the integrated power function through a two-sided sandwich bound (Section~\ref{sec:sandwich}). Its lower side carries a Davis--Kahan-type dependence~\cite{DavisKahan1970} on the spectral gap, which explains why a small Nystr\"om trace error certifies subspace alignment only when the spectrum has a clear gap at $s$.

The second group of results explains why the minimizers are evenly spread (Section~\ref{sec:evenness}),
which we observed in every experiment we ran on $\Sph^1$ and $\Sph^2$. For $m=s$, we write $F$ through the tangents of the principal angles between the landmark space and $V_s$ (Proposition~\ref{thm:factor}, a direct application of \cite{ZhuKnyazev2013}). This expresses $F$ through the aliasing of the tail eigenfunctions onto $V_s$ and yields an upper bound that is minimized exactly by \emph{tight sampling sets}, meaning configurations whose eigenfunction sampling matrix is a scaled isometry (Theorem~\ref{thm:tight}). In a homogeneous setting, such sets are precisely the minimal equal-weight cubature rules for products of leading eigenfunctions (Theorem~\ref{thm:tightchar}), a homogeneous form of the classical correspondence between minimal cubature and common zeros of orthogonal polynomials \cite{Mysovskikh1981,CoolsMysovskikhSchmid2001}. Equivalently, the sampled eigenfunction vectors form a unit-norm tight frame, the configuration minimizes the frame potential, and for $m=s$ it is a set of Fekete points for $V_s$ attaining the Hadamard bound (Proposition~\ref{prop:frames}). The frame description also extends to every $m\ge s$. Tight sets are the equispaced points on $\Sph^1$ and the regular tetrahedron on $\Sph^2$ for $s=4$ (Corollaries~\ref{cor:tight_circle} and~\ref{cor:tight_sphere}). For larger $s$ on $\Sph^2$, tightness with $m=s$ is impossible \cite{BannaiDamerell1979}. A complementary symmetry principle, a special case of Michel's theorem on critical orbits \cite{Michel1971}, shows that any orbit of a finite isometry group whose stabilizers fix no tangent vector is a critical point of $F$ (Theorem~\ref{thm:symmetry}). This accounts for the equispaced circle configurations and the Platonic solids. Criticality does not imply minimality, however, as the cube for $s=4$ and the dodecahedron for $s=9$ show, even though both are tight frames for those values of $s$. We also give a closed-form aliasing formula on the circle, which gives the value of $F$ at equispaced landmarks for every $m\ge s$ and decays super-exponentially in $m$ (Theorem~\ref{thm:alias}).

Finally, we study the numerical behavior of the problem and test the theory experimentally. The Gaussian bandwidth must lie in a window proportional to the landmark spacing $h$. Below roughly $0.15\,h$ the position-dependent part of $F$ becomes exponentially small and falls below the resolution of the optimizer. At the upper end two mechanisms intervene. Once the attainable minimum of $F$ drops below about $10^{-9}$, it is poorly resolved by the subtraction $s-\tr(\cdot)$ in double precision and absolute stopping tolerances trigger early. At larger bandwidths, beyond roughly $2.4\,h$ on the circle and somewhat less on the sphere, the landmark Gram matrix exceeds a condition number of $10^{12}$. The resulting window shrinks as landmarks are added (Section~\ref{sec:conditioning}). Within this window, every random start on the circle converges to equispaced landmarks, and on the sphere the optimizer recovers the regular tetrahedron, octahedron, and icosahedron (Sections~\ref{sec:circle}--\ref{sec:sphere}). The optimized sets also have markedly smaller mesh ratios than Fibonacci lattices, and the integrated power function at the optimized landmarks sits near the lower end of the sandwich bound.

\section{Background and Related work}
The Nystr\"om method is named after Nystr\"om's quadrature treatment of integral
equations \cite{Nystrom1930,Baker1977}, and entered machine learning as a way
to speed up kernel machines \cite{WilliamsSeeger2001,SmolaScholkopf2000,FineScheinberg2001}.
Many subsequent analyses bound the error of $K-\widehat K$ in the trace, spectral
or Frobenius norm and relate it to the choice of landmarks
\cite{DrineasMahoney2005,ZhangTsangKwok2008,KumarMohriTalwalkar2012,GittensMahoney2016}. Choosing the best $m$ landmarks is a combinatorial problem
\cite{OglicGartner2017}, and a range of strategies have been proposed. These
include uniform sampling \cite{WilliamsSeeger2001}, leverage-score and ridge
leverage-score sampling \cite{DrineasMahoney2005,AlaouiMahoney2015,MuscoMusco2017},
$k$-means centers \cite{ZhangTsangKwok2008,OglicGartner2017}, determinantal
sampling \cite{BelabbasWolfe2009}, regularized Christoffel functions
\cite{FanuelSchreursSuykens2022}, and greedy pivoted Cholesky
\cite{SmolaScholkopf2000,FineScheinberg2001,HarbrechtPetersSchneider2012}.
The guarantees for these schemes are stated for the residual $K-\widehat K$,
not for the alignment of $\Span\{k(\cdot,z_j)\}$ with an eigenspace of $K$.

The trace error is exactly $\int P_Z^2\dd$ (Section~\ref{sec:intpower}), which ties the Nystr\"om literature to scattered-data approximation, where
the power function controls the pointwise error of radial basis function
(RBF) interpolation \cite{Schaback1995,SchabackWendland2006,Wendland2005}. In that
setting, point sets that are near-optimal for the power function are
quasi-uniform \cite{DeMarchiSchabackWendland2005}, and $P$-greedy selection
attains optimal convergence rates \cite{SantinHaasdonk2017}. Fekete points, which maximize the
determinant of a sampling matrix, offer a classical route to well-conditioned
interpolation along similar lines \cite{BosEtAl2010}. Our objective differs
from all of these in that it is a spectral subspace criterion rather than a
residual norm. The sandwich bound (Theorem~\ref{thm:sandwich}) makes the
relation between the two quantitative, and Theorem~\ref{thm:tight} shows that
spectral alignment rewards a property closely related to the Fekete idea,
namely an isometric eigenfunction sampling matrix.

Because evenly spaced landmarks emerge from the optimization, our sphere
experiments connect to the standard notions of uniformity on $\Sph^2$.
Spherical designs \cite{DelsarteGoethalsSeidel1977, BannaiDamerell1979} and minimal-energy or best-packing configurations \cite{SaffKuijlaars1997}, the latter going back to the Tammes problem \cite{Tammes1930,FejesToth1943,SchutteVanDerWaerden1951,Danzer1986}, are the usual notions of evenly spaced points on $\Sph^2$. The Platonic solids are the canonical small examples, and our sphere experiments recover several of them.

Finally, the practical limits we identify for Gaussian kernels belong to a
well-studied trade-off between accuracy and conditioning for flat RBFs, which
Schaback's uncertainty principle describes \cite{Schaback1995}. Stable
evaluation schemes such as RBF-QR were developed to circumvent it
\cite{FornbergPiret2008,FornbergLarssonFlyer2011}, and they form the
background for Section~\ref{sec:conditioning}.

\subsection{Setting and conventions}\label{sec:setting}

Let $\Omega$ be a compact metric space, $\mu$ a finite Borel measure on $\Omega$, and $K$ a continuous, symmetric, strictly positive-definite kernel on $\Omega$ with reproducing kernel Hilbert space (RKHS) $\Hcal$. Write $K_x=K(\cdot,x)$, so that $\langle f,K_x\rangle_{\Hcal}=f(x)$ for all $f\in\Hcal$ (the reproducing property). The integral operator $(T_Kf)(x)=\int K(x,y)f(y)\dd(y)$ is compact, self-adjoint and positive on $L^2(\mu)$. We write its positive eigenvalues as $\lambda_1\ge\lambda_2\ge\dots>0$, with $L^2(\mu)$-orthonormal eigenfunctions $\varphi_i$. There are infinitely many of them when $\mu$ has infinite support, as for the uniform measures on the circle and the sphere, and at most $N$ of them when $\mu$ is an empirical measure on $N$ points. Each $\varphi_i$ has the continuous representative $\varphi_i=\lambda_i^{-1}T_K\varphi_i$, which is defined on all of $\Omega$, and Mercer's theorem gives
\[K(x,y)=\sum_{i}\lambda_i\,\varphi_i(x)\varphi_i(y),\qquad x,y\in\operatorname{supp}\mu .
\] The functions
\[
\psi_i:=\sqrt{\lambda_i}\,\varphi_i,\qquad \langle\psi_i,\psi_k\rangle_{\Hcal}=\delta_{ik},
\]
are orthonormal in $\Hcal$, and the reproducing property gives $\langle\psi_i,K_x\rangle_{\Hcal}=\psi_i(x)$ for every $x\in\Omega$.

\vspace{3mm}

The $\psi_i$ need not span $\Hcal$. Their closed span is the closure of the range of the covariance operator $C$ of Section~\ref{sec:intpower}, and its orthogonal complement is
\[\Hcal_0:=\{f\in\Hcal:\ f=0\ \ \mu\text{-a.e.}\}=\ker C .
\]
When $\operatorname{supp}\mu=\Omega$, continuity forces $\Hcal_0=\{0\}$, so $\{\psi_i\}$ is an orthonormal basis of $\Hcal$. For an empirical measure $\mu_N$, by contrast, $\Hcal_0$ is infinite-dimensional, and the translate $K_z$ of a point $z$ off the data has a nonzero component in it. To treat both cases at once we complete $\{\psi_i\}$ to an orthonormal basis of $\Hcal$ by appending an orthonormal basis of $\Hcal_0$, and we assign the eigenvalue $\lambda_i=0$ to each appended element. These are placed after all positive eigenvalues, so the ordering $\lambda_1\ge\lambda_2\ge\dots\ge0$ is preserved. With this convention $\{\psi_i\}$ is always an orthonormal basis of $\Hcal$ that diagonalizes $C$, the expansion
\[K_x=\sum_i\psi_i(x)\,\psi_i
\] holds for every $x\in\Omega$, and every sum over $i$ below runs over the completed basis. Throughout we assume $\lambda_s>0$, so that the leading eigenspace $V_s$ defined below consists of genuine eigenfunctions of $T_K$.

\vspace{3mm}

\noindent For the landmarks $Z=\{z_1,\dots,z_m\}$ with $K_{ZZ}=[K(z_j,z_l)]$ invertible, define the two subspaces:
\begin{itemize}[leftmargin=*]
\item $U_Z=\Span\{K_{z_1},\dots,K_{z_m}\}$, with $\Hcal$-orthogonal projector $\Pi_Z$.
\item $V_s=\Span\{\psi_1,\dots,\psi_s\}=\Span\{\varphi_1,\dots,\varphi_s\}$, with $\Hcal$-orthogonal projector $Q_s$.
\end{itemize}

\subsection{The Projector and the Power function}\label{sec:projector}
For landmarks $Z=\{z_1,\dots,z_m\}$ with $K_{ZZ}$ invertible, the orthogonal projector onto the span of the kernel translates is the operator that interpolates a function at the landmarks using the translates $K_{z_j}$, namely
\begin{equation}\label{eq:Pi}
\Pi_Zf=\sum_{j,l=1}^m K_{z_j}\,(K_{ZZ}^{-1})_{jl}\,f(z_l)=k_Z(\cdot)^\top K_{ZZ}^{-1}f(Z),
\end{equation}
where $k_Z(\cdot)=(K(\cdot,z_j))_{j=1}^m$ and $f(Z)=(f(z_1),\dots,f(z_m))^\top$. To see that this is the orthogonal projector onto $U_Z=\Span\{K_{z_1},\dots,K_{z_m}\}$, we check three properties. Its range lies in $U_Z$ by construction. It fixes $U_Z$, because applying it to $f=K_{z_k}$ gives $f(z_l)=K(z_l,z_k)=(K_{ZZ})_{lk}$ and the coefficient sum collapses to $\sum_j K_{z_j}\sum_l (K_{ZZ}^{-1})_{jl}(K_{ZZ})_{lk}=K_{z_k}$. And it is self-adjoint, since the reproducing property gives $\langle g,\Pi_Zf\rangle_{\Hcal}=g(Z)^\top K_{ZZ}^{-1}f(Z)$, which is symmetric in $f$ and $g$.

The same computation yields a useful formula for the size of a projected function. Because $\Pi_Z$ is idempotent and self-adjoint, $\|\Pi_Zf\|_{\Hcal}^2=\langle f,\Pi_Zf\rangle_{\Hcal}$, and applying $\langle f,K_{z_j}\rangle_{\Hcal}=f(z_j)$ to \eqref{eq:Pi} turns this into
\begin{equation}\label{eq:quadform}
\|\Pi_Zf\|_{\Hcal}^2=f(Z)^\top K_{ZZ}^{-1}f(Z).
\end{equation}
In other words, the portion of $f$ captured by the landmark space depends only on the values of $f$ at the landmarks, weighted by the inverse Gram matrix. This identity does most of the work in what follows, because it lets us express every projected norm through finite-dimensional matrices.

The power function arises when we apply this to the feature vectors themselves. We define $P_Z(x)^2:=\|K_x-\Pi_ZK_x\|_{\Hcal}^2$, which measures how far the translate at $x$ lies from the landmark space. Pythagoras gives $\|u-\Pi u\|^2=\|u\|^2-\|\Pi u\|^2$, and combining this with $\|K_x\|_{\Hcal}^2=K(x,x)$ and \eqref{eq:quadform} applied to $f=K_x$ (for which $f(Z)=k_Z(x)$) produces
\begin{equation}\label{eq:power}
P_Z(x)^2=K(x,x)-k_Z(x)^\top K_{ZZ}^{-1}k_Z(x).
\end{equation}
This is the classical power function formula \cite{Wendland2005}. It is also the diagonal of the Nystr\"om residual $K-\widehat K$, where $\widehat K(x,y)=k_Z(x)^\top K_{ZZ}^{-1}k_Z(y)$, so the same object serves as the pointwise error of kernel interpolation and as the pointwise error of the Nystr\"om approximation. The remainder of the paper compares two ways of summarizing the quality of $Z$. One integrates $P_Z^2$ over the domain (Section~\ref{sec:intpower}), and the other measures the alignment of $U_Z$ with $V_s$ directly through the spectral alignment $F(Z)$ (Section~\ref{sec:discrepancy}).

\section{The Spectral Alignment functional}\label{sec:discrepancy}

Our goal is to choose landmarks so that $U_Z$ contains $V_s$ as closely as possible, so we first need a measure of how far $V_s$ lies from the landmark space. When $m=s$ the two subspaces have the same dimension and the Hilbert--Schmidt distance between their projectors is a natural choice. When $m>s$ that distance can never vanish, because the subspaces differ in dimension, so we measure instead how much of $V_s$ is missing from $U_Z$. Expanding the square and using $\Pi_Z^2=\Pi_Z$, $Q_s^2=Q_s$, and cyclicity of the trace gives
\[
\|\Pi_Z-Q_s\|_{HS}^2=\tr\big((\Pi_Z-Q_s)^2\big)=\tr\Pi_Z+\tr Q_s-2\tr(\Pi_ZQ_s)=m+s-2\tr(\Pi_ZQ_s),
\]
since $\tr\Pi_Z=m$ and $\tr Q_s=s$. We define the spectral alignment functional to be
\[
F(Z):=s-\tr(\Pi_ZQ_s).
\]
When $m> s$ this quantity is related to the projector distance by $\tfrac12\|\Pi_Z-Q_s\|_{HS}^2=F(Z)+\tfrac{m-s}{2}$. The two coincide when $m=s$, and for $m>s$ the alignment $F$ removes the constant offset that the dimension mismatch forces on the Hilbert--Schmidt distance.

The geometric meaning follows from the principal angles. The operator $Q_s\Pi_ZQ_s$ restricted to $V_s$ has eigenvalues $\cos^2\theta_j$ for $j=1,\dots,s$, where $\theta_1,\dots,\theta_s$ are the principal angles between $V_s$ and $U_Z$. Hence $\tr(\Pi_ZQ_s)=\sum_j\cos^2\theta_j$ and
\[
F(Z)=\sum_{j=1}^s\sin^2\theta_j\in[0,s].
\]
This is the standard chordal Grassmannian distance, squared, between $V_s$ and the subspace of $U_Z$ closest to it. It vanishes exactly when $V_s\subseteq U_Z$, which for $m=s$ means $U_Z=V_s$.

To understand where the spectral alignment comes from, it helps to track how much of each eigenvector lies inside the landmark space. We define the overlap coefficients $a_i:=\|\Pi_Z\psi_i\|_{\Hcal}^2\in[0,1]$. Applying \eqref{eq:quadform} with $f=\psi_i$ and using $\psi_i=\sqrt{\lambda_i}\varphi_i$ gives
\begin{equation}\label{eq:overlap}
a_i=\psi_i(Z)^\top K_{ZZ}^{-1}\psi_i(Z)=\tr\big(K_{ZZ}^{-1}\,\lambda_i\varphi_i(Z)\varphi_i(Z)^\top\big).
\end{equation}
We can think of $V_s$ as the target subspace and $U_Z$ as the one we are able to build from landmarks. Each $\psi_i$ with $i\le s$ should lie entirely inside $U_Z$, and the amount by which it fails is $1-a_i$. Computing $\tr(\Pi_ZQ_s)$ in the orthonormal basis $\{\psi_i\}$, and noting that $Q_s\psi_i=\psi_i$ for $i\le s$ and $Q_s\psi_i=0$ otherwise, we find $\tr(\Pi_ZQ_s)=\sum_{i\ge1}\langle\psi_i,\Pi_ZQ_s\psi_i\rangle = \sum_{i\le s}\langle\psi_i,\Pi_Z\psi_i\rangle=\sum_{i\le s}a_i$. The spectral alignment is therefore the total leakage from $V_s$,
\[
F(Z)=s-\sum_{i\le s}a_i=\sum_{i\le s}(1-a_i).
\]

The same coefficients also reveal a second interpretation. Taking the trace of $\Pi_Z$ in the basis $\{\psi_i\}$ gives $\sum_{i\ge1}a_i=\tr\Pi_Z=m$, and subtracting the first leaves
\begin{equation}\label{eq:massconservation}
\sum_{i>s}a_i=m-s+F(Z).
\end{equation}
The tail $V_s^\perp$ therefore receives at least $m-s$ units of mass, since the landmark space has $m$ dimensions and only $s$ of them are needed to cover $V_s$. The alignment $F(Z)$ is exactly the amount by which the tail mass exceeds this unavoidable minimum, and it is simultaneously the mass of $V_s$ that leaks out of $U_Z$. When $m=s$ the minimum is zero and the two views coincide, so $F(Z)=\sum_{i>s}a_i$, that is, the mass that $V_s$ loses is exactly the mass that the tail gains. The first view asks whether the landmarks can reproduce the leading eigenfunctions, while the second asks how much tail mass the landmark space holds beyond the unavoidable $m-s$.

The overlap coefficients can be summed in closed form without computing any eigenfunctions explicitly. We split the kernel into its truncated and tail parts,
\[
K^{(s)}=\sum_{i\le s}\lambda_i\varphi_i\otimes\varphi_i,\qquad T^{(s)}=K-K^{(s)}=\sum_{i>s}\lambda_i\varphi_i\otimes\varphi_i,
\]
and write $K^{(s)}_{ZZ}$ and $T^{(s)}_{ZZ}$ for their $m\times m$ landmark matrices, so that $K_{ZZ}=K^{(s)}_{ZZ}+T^{(s)}_{ZZ}$. Equivalently, $K^{(s)}=\sum_{i\le s}\psi_i\otimes\psi_i$ and $T^{(s)}=\sum_{i>s}\psi_i\otimes\psi_i$ are the reproducing kernels of $V_s$ and of $V_s^\perp$ in $\Hcal$, so both are positive semidefinite. Summing \eqref{eq:overlap} over $i\le s$ gives $\sum_{i\le s}a_i=\tr(K_{ZZ}^{-1}K^{(s)}_{ZZ})$, and therefore
\begin{equation}\label{eq:F}
\boxed{\ F(Z)=s-\tr\big(K_{ZZ}^{-1}K^{(s)}_{ZZ}\big)\ }
\end{equation}
Since $\tr(K_{ZZ}^{-1}K_{ZZ})=m$, the same quantity can be written in terms of the tail kernel as
\begin{equation}\label{eq:closedform}
F(Z)=\tr\big(K_{ZZ}^{-1}T^{(s)}_{ZZ}\big)-(m-s).
\end{equation}
When $m=s$ the correction vanishes and $F(Z)=\tr(K_{ZZ}^{-1}T^{(s)}_{ZZ})$, whose nonnegativity is immediate because $K_{ZZ}^{-1}\succ0$ and $T^{(s)}_{ZZ}\succeq0$. A quick sanity check at $s=m=1$ confirms the formulas. For a single landmark $z$ we have $\cos^2\theta=\langle\psi_1,K_z\rangle_{\Hcal}^2/\|K_z\|_{\Hcal}^2=\lambda_1\varphi_1(z)^2/K(z,z)$, so $F=1-\lambda_1\varphi_1(z)^2/K(z,z)=T^{(1)}(z,z)/K(z,z)$, which is exactly what \eqref{eq:closedform} gives when $K_{ZZ}=K(z,z)$.

\begin{remark}[Basis invariance]\label{rem:basis}
Only $K^{(s)}_{ZZ}=\Psi\Psi^\top$ enters \eqref{eq:F}, where $\Psi_{ji}=\psi_i(z_j)$. Hence $F$ does not depend on which orthonormal basis of $V_s$ is used. This matters in two places. Inside a degenerate eigenspace any rotation of the basis is allowed, and on the circle and the sphere (Appendix~\ref{sec:homog}) the kernel $K^{(s)}$ can be written down without ever forming a basis.
\end{remark}

The formulas also tell us when the alignment can vanish. For $m=s$ we have $F(Z)=0$ if and only if $T^{(s)}_{ZZ}=0$, and since $T^{(s)}_{ZZ}$ is positive semidefinite this forces its diagonal to vanish. That requires $\psi_i(z_j)=0$ for every $i>s$ and every $j$, so the landmarks would have to be common zeros of all the tail eigenfunctions. For a kernel with an infinite spectrum this common-zero condition is not expected to hold in general, so $\min_Z F(Z)>0$ generically and the problem is one of minimization rather than exact solution. For $m>s$ the condition for $F(Z)=0$ is instead $V_s\subseteq U_Z$, which requires every leading eigenfunction to be an exact combination of $m$ kernel translates. This is also generically impossible for an infinite spectrum, but unlike the case $m=s$ it can be approached closely. The aliasing formula of Theorem~\ref{thm:alias} shows that on the circle $F$ decays super-exponentially in $m$ at equispaced landmarks.

The closed form is also practical for data. Let $\mu$ be the empirical measure on $N$ points, with kernel matrix $K_N$ and rank-$s$ truncation $K_N^{(s)}=U_s\Sigma_sU_s^\top$. For landmark set $Z$ with $|Z|=m\ge s$, the spectral alignment is $F(Z)=s-\tr\big((K_N)_{ZZ}^{-1}(K_N^{(s)})_{ZZ}\big)$. The $N^{-1}$ normalizations of $\lambda_i$ and $\varphi_i$ cancel in this expression, so no rescaling is needed. Evaluating it requires only the top $s$ eigenpairs of $K_N$ and the two $m\times m$ submatrices for each candidate set. This form assumes that the landmarks are data points. Landmarks off the data are handled in Section~\ref{sec:algorithm} through the Nystr\"om extension, and there the completed basis of Section~\ref{sec:setting} is essential, because the translates $K_{z_j}$ then have components in $\Hcal_0$.

\subsection{The integrated power function as a weighted distance}\label{sec:intpower}

The power function \eqref{eq:power} is a pointwise quantity, and to compare it with the alignment functional we need a single number that summarizes it over the domain. The natural choice is its integral against $\mu$, which equals the Nystr\"om trace error. To express this integral in the same Hilbert space language as the previous section, we introduce the covariance operator $C:\Hcal\to\Hcal$,
\[
Cf=\int f(x)\,K_x\dd(x)=\int\langle f,K_x\rangle_{\Hcal}\,K_x\dd(x),
\qquad\text{that is,}\qquad
C=\int K_x\otimes K_x\dd(x).
\]
Each integrand $K_x\otimes K_x$ is a rank-one operator, so $C$ is the average of these operators and plays the role of the second-moment operator of the feature vectors $K_x\in\Hcal$. Its role is the same as that of $\frac{1}{N}\sum_n x_nx_n^\top$ for vectors in $\R^d$. The integral is a Bochner integral in $\Hcal$, and it exists because $\|K_x\|_{\Hcal}^2=K(x,x)$ is bounded for a continuous kernel on a compact space and $\mu$ is finite.

Evaluating $Cf$ at a point $y$ gives $(Cf)(y)=\int K(y,x)f(x)\dd(x)$, which is the integral operator $T_K$ whose eigenpairs define the Mercer expansion, so $C$ is simply $T_K$ viewed as an operator on $\Hcal$. Its matrix in the completed orthonormal basis $\{\psi_i\}$ is diagonal, since
\[
\langle\psi_i,C\psi_k\rangle_{\Hcal}=\int\psi_i(x)\psi_k(x)\dd(x)=\sqrt{\lambda_i\lambda_k}\,\delta_{ik}=\lambda_i\delta_{ik}.
\]
Therefore $C=\sum_i\lambda_i\,\psi_i\otimes\psi_i$ and $\tr C=\sum_i\lambda_i=\int K(x,x)\dd(x)$.

With $C$ in hand, the integrated power function becomes a trace. Writing $P_Z(x)^2=\langle K_x,(I-\Pi_Z)K_x\rangle_{\Hcal}=\tr\big((I-\Pi_Z)\,K_x\otimes K_x\big)$ and integrating over $x$, we obtain
\[
\int P_Z^2\dd=\tr\big((I-\Pi_Z)C\big)=\big\|(I-\Pi_Z)C^{1/2}\big\|_{HS}^2.
\]
This identity shows that the integrated power function is a Hilbert--Schmidt distance, but between the operator $C^{1/2}$ and its projection onto $U_Z$. Since $C^{1/2}$ is not a projector (it has full rank and weights each eigen-direction by $\sqrt{\lambda_i}$), this is not a distance between $U_Z$ and any fixed $s$-dimensional subspace. Our spectral alignment $F(Z)=\|(I-\Pi_Z)Q_s\|_{HS}^2$ has the same form, but with the spectral projector $Q_s$ in place of $C^{1/2}$. The difference becomes explicit when we evaluate the trace in the eigenbasis. Since $\langle\psi_i,(I-\Pi_Z)\psi_i\rangle_{\Hcal}=1-a_i$, we find
\begin{equation}\label{eq:intpower_eig}
\int P_Z^2\dd=\sum_{i\ge1}\lambda_i\,(1-a_i).
\end{equation}
The integrated power function thus weights the unexplained mass of each eigenvector by its eigenvalue, whereas the alignment counts every leading eigenvector equally.

Equation \eqref{eq:intpower_eig} also gives a lower bound that no choice of landmarks can beat. The coefficients satisfy $0\le a_i\le1$ and $\sum_ia_i=m$, so maximizing $\sum_i\lambda_ia_i$ under these constraints places $a_i=1$ on the $m$ largest eigenvalues, which gives $\sum_{i\le m}\lambda_i$ by the Ky Fan maximum principle \cite{KyFan1949}. Consequently
\begin{equation}\label{eq:kyfan}
\int P_Z^2\dd\ \ge\ \sum_{i>m}\lambda_i,
\end{equation}
with equality if and only if $U_Z=\Span\{\psi_1,\dots,\psi_m\}$, provided $\lambda_m>\lambda_{m+1}$. The right-hand side is the error of the best rank-$m$ approximation, so \eqref{eq:kyfan} says that $m$ landmarks cannot do better than the truncated Mercer expansion of the same rank.

Because our target is the rank-$s$ eigenspace, it is natural to measure the quality of a landmark set against the rank-$s$ benchmark and define the excess $E(Z):=\int P_Z^2\dd-\sum_{i>s}\lambda_i$. When $m=s$ the bound \eqref{eq:kyfan} makes $E(Z)\ge0$. When $m>s$ the excess can be negative, since extra landmarks can capture part of the tail, but \eqref{eq:kyfan} still gives $E(Z)\ge-\sum_{s<i\le m}\lambda_i$. Splitting \eqref{eq:intpower_eig} at index $s$ gives, for every $m\ge s$,
\begin{equation}\label{eq:excess}
E(Z)=\sum_{i\le s}\lambda_i(1-a_i)-\sum_{i>s}\lambda_ia_i=\tr\big((Q_s-\Pi_Z)C\big).
\end{equation}
Comparing this with the alignment functional shows that the same leakage numbers appear in both quantities. The spectral alignment is $F(Z)=\sum_{i\le s}(1-a_i)$, and the tail coefficients satisfy $\sum_{i>s}a_i=m-s+F(Z)$ by \eqref{eq:massconservation}. The excess, however, weights these numbers by the eigenvalues and gives the two groups opposite signs. Leakage out of the leading eigenspace is penalized in proportion to the large eigenvalues $\lambda_i$ with $i\le s$, while mass of the tail captured by $U_Z$ offsets it in proportion to the small eigenvalues with $i>s$. This opposition is what prevents a small excess from forcing a small alignment when the spectrum has no gap, and the extra tail mass $m-s$ that the landmark space must hold when $m>s$ enters the bound of the next section.

\subsection{The Sandwich Bound}\label{sec:sandwich}

The spectral alignment $F$ and the excess $E$ are built from the same leakage numbers, so it is natural to ask how tightly they constrain one another. To compare them, write $b_i=1-a_i\ge0$ for $i\le s$. The definition of $F$ gives $\sum_{i\le s}b_i=F(Z)$, and \eqref{eq:massconservation} gives $\sum_{i>s}a_i=m-s+F(Z)$. In this notation \eqref{eq:excess} reads $E=\sum_{i\le s}\lambda_ib_i-\sum_{i>s}\lambda_ia_i$. The eigenvalues are ordered, so $\lambda_i\ge\lambda_s$ for $i\le s$, and the first sum is at least $\lambda_sF$. For the second sum we also use $a_i\le1$. The tail mass $\sum_{i>s}a_i=m-s+F$ splits into the part on the indices $s<i\le m$ and the part on $i>m$, where every eigenvalue is at most $\lambda_{m+1}$. Hence\[
\sum_{i>s}\lambda_ia_i\ \le\ \sum_{s<i\le m}(\lambda_i-\lambda_{m+1})\,a_i+\lambda_{m+1}(m-s+F)\ \le\ \sum_{s<i\le m}\lambda_i+\lambda_{m+1}F,
\]
where the last step uses $a_i\le1$ and $\lambda_i\ge\lambda_{m+1}$ for $i\le m$. The extremal case fills the directions $\psi_{s+1},\dots,\psi_m$ completely and places the leaked mass $F$ on $\psi_{m+1}$. Combining the two sums yields the lower estimate
 \[
E\ \ge\ \lambda_sF-\sum_{s<i\le m}\lambda_i-\lambda_{m+1}F=(\lambda_s-\lambda_{m+1})\,F-\sum_{s<i\le m}\lambda_i .
\] 
For the upper estimate, we may discard the nonpositive second sum and use $\lambda_i\le\lambda_1$ in the first, which gives $E\le\sum_{i\le s}\lambda_ib_i\le\lambda_1F$. Together these two estimates constitute the following result.

\begin{theorem}[Sandwich bound]\label{thm:sandwich}
For $m\ge s$ landmarks with $K_{ZZ}$ invertible,
\[
\boxed{(\lambda_s-\lambda_{m+1})\,F(Z)-\sum_{s<i\le m}\lambda_i\ \le\ \int P_Z^2\dd-\sum_{i>s}\lambda_i\ \le\ \lambda_1\,F(Z).}
\]
In particular, for $m=s$ the bound reads $(\lambda_s-\lambda_{s+1})\,F(Z)\le E(Z)\le\lambda_1F(Z)$. Since $\sum_{s<i\le m}\lambda_i\le(m-s)\lambda_{s+1}$ and $\lambda_{m+1}\le\lambda_{s+1}$, the lower bound also implies the simpler estimate $E(Z)\ge(\lambda_s-\lambda_{s+1})F(Z)-\lambda_{s+1}(m-s)$. \end{theorem}

The upper inequality says that driving $F$ to zero always keeps the trace error within $\lambda_1F$ of the best rank-$s$ approximation, which is the sense in which small leakage implies a near-optimal Nystr\"om approximation. The lower inequality runs in the opposite direction and is informative only when the spectrum has a gap. For $m=s$, a small Nystr\"om trace excess implies $F\le E/(\lambda_s-\lambda_{s+1})$, so alignment follows from near-optimal trace error only if $\lambda_s-\lambda_{s+1}$ is not small. Without a gap, $U_Z$ can be far from $V_s$ while $\int P_Z^2\dd$ is nearly optimal. This is the Davis--Kahan phenomenon \cite{DavisKahan1970}, in which the stability of an eigenspace under perturbation scales inversely with the gap to the next eigenvalue.

For $m>s$ the lower bound carries the additional correction $\sum_{s<i\le m}\lambda_i$. It arises because the landmark space has $m-s$ dimensions beyond those needed to cover $V_s$, and these can absorb the eigen-directions $\psi_{s+1},\dots,\psi_m$, each of which offsets the excess by $\lambda_i$. At $F=0$ the bound reduces to $E\ge-\sum_{s<i\le m}\lambda_i$, which is exactly the Ky Fan bound \eqref{eq:kyfan}.

A concrete case shows that, for $m=s$, the constant $\lambda_s-\lambda_{s+1}$ in the lower bound cannot be improved by any bound that holds for all $s$-dimensional subspaces of $\Hcal$. Suppose $U_Z$ differs from $V_s$ by a single rotation through an angle $\theta$ in the plane spanned by $\psi_s$ and $\psi_{s+1}$, so that $U_Z$ contains $\cos\theta\,\psi_s+\sin\theta\,\psi_{s+1}$ in place of $\psi_s$ and the rest of $V_s$ is unchanged. The overlap coefficients then satisfy $a_s=\cos^2\theta$ and $a_{s+1}=\sin^2\theta$, with all other $a_i$ unchanged. It follows that $F=1-a_s=\sin^2\theta$ and $E=\lambda_s\sin^2\theta-\lambda_{s+1}\sin^2\theta=(\lambda_s-\lambda_{s+1})\sin^2\theta$, so $E=(\lambda_s-\lambda_{s+1})F$ holds with equality. At $\theta=90^\circ$ the subspace $U_Z$ has swapped $\psi_s$ for $\psi_{s+1}$ and $F=1$, as large as it can be for a one-direction swap. If the gap $\lambda_s-\lambda_{s+1}=\varepsilon$ is tiny, however, the excess is only $E=\varepsilon$. The trace error is then almost optimal even though the subspace has lost an entire leading direction, so near-optimality of $\int P_Z^2\dd$ says almost nothing about alignment with $V_s$.


It is worth separating this near-degenerate situation from the exactly degenerate one. If $\lambda_s=\lambda_{s+1}$, the subspace $V_s$ is not uniquely defined and $F$ depends on an arbitrary choice. The swapped subspace in the example above is then an equally valid top-$s$ eigenspace, and $E=0$ correctly reports that it is optimal. The difficulty lies in the near-degenerate case, where $V_s$ is well defined but any rotation within the cluster costs only $\varepsilon\sin^2\theta$ in excess. The practical conclusion is that small $F$ always forces a small $E$, but a small $E$ forces a small $F$ only when there is a gap between $\lambda_s$ and $\lambda_{s+1}$. For landmark determination this means that minimizing $\int P_Z^2\dd$ is a legitimate proxy for subspace alignment only when the spectrum has a clear gap at $s$.

The same leakage numbers also give direct estimates for the integrated power function itself, which are useful when comparing with the rank-$s$ benchmark. They follow from \eqref{eq:intpower_eig} directly rather than from Theorem~\ref{thm:sandwich}.

\begin{corollary}[Bounds for $m\ge s$]\label{cor:sgeqm}
For $m\ge s$ landmarks with $K_{ZZ}$ invertible,
\[
\lambda_s\,F(Z)\ \le\ \int P_Z^2\dd\ \le\ \sum_{i>s}\lambda_i+\lambda_1F(Z),\qquad \int P_Z^2\dd\ \ge\ \sum_{i>m}\lambda_i .
\]
\end{corollary}


\section{Evenly Spaced Minimizers}\label{sec:evenness}

In every experiment we ran on $\Sph^1$ and $\Sph^2$ the minimizers of $F$ turned out to be evenly spread, even though nothing in the objective asks for this. This section explains the phenomenon by connecting $F$ to three classical families of evenly spread point sets: minimal positive cubature rules and tight spherical designs, unit-norm tight frames, and Fekete points. Most of the ingredients are known, and we cite them as such. What is new is the link between these objects and the spectral alignment, together with the quantitative bound of Theorem~\ref{thm:tight} and the closed form of Theorem~\ref{thm:alias}.

The argument has three steps, all valid on any compact homogeneous space with an invariant kernel. First, for $m=s$, we write $F$ through the tangents of the principal angles between $U_Z$ and $V_s$. This yields an upper bound controlled by the conditioning of the eigenfunction sampling matrix, and the bound is minimized exactly by \emph{tight sampling sets}. Second, we identify tight sampling sets with minimal equal-weight cubature rules, with unit-norm tight frames, and with Fekete points. The tight-frame description extends to every $m\ge s$. Third, we recall that sufficiently symmetric configurations are critical points of every invariant objective, and hence of $F$. Facts specific to the circle and the sphere are collected in Appendix~\ref{app:circlesphere}.

\subsection{The Homogeneous Setting}\label{sec:homog_general}

We say that $(\Omega,\mu,K,G)$ is a \emph{homogeneous setting} if $\Omega$ is a compact Hausdorff space (for instance a compact Riemannian manifold), $G$ is a group of homeomorphisms of $\Omega$ that acts transitively, $\mu$ is a $G$-invariant Borel probability measure, and $K$ is a continuous, symmetric, strictly positive-definite kernel with $K(gx,gy)=K(x,y)$ for all $g\in G$ and $x,y\in\Omega$.

Examples include $\Sph^d$ and $\mathbb{RP}^d$ with $G$ the orthogonal group, flat tori with $G$ the translations, compact Lie groups acting on themselves by left translation, and Grassmann and Stiefel manifolds. Two families of kernels satisfy the hypotheses on every such space. The first is the heat kernel of an invariant Riemannian metric, and the second is the ambient Gaussian $K(x,y)=\exp(-\|x-y\|^2/2\ell^2)$ restricted to an orbit $\Omega=G\cdot x_0\subset\R^n$ of an orthogonal representation of a compact group. The latter is invariant because $G$ acts by orthogonal maps, and it is strictly positive definite because the Euclidean Gaussian is. The Gaussian kernels on $\Sph^1$ and $\Sph^2$ used in our experiments are of the second type, and their exact Mercer expansions are recorded in Appendix~\ref{sec:homog}. By contrast, the Gaussian of the \emph{geodesic} distance is positive definite for all bandwidths only when the manifold is isometric to a Euclidean space \cite{FeragenLauzeHauberg2015}. This excludes not only curved spaces but also the circle and flat tori, which are flat but compact, so the geodesic Gaussian should not be used on any of the compact homogeneous spaces considered here.

For $g\in G$ define $(\pi_gf)(x)=f(g^{-1}x)$. Invariance of $\mu$ makes $\pi_g$ an orthogonal operator on $L^2(\mu)$, and invariance of $K$ makes it commute with the integral operator $T_K$. Hence every eigenspace $\mathcal E_\lambda$ of $T_K$ is $\pi_g$-invariant. Eigenfunctions with $\lambda>0$ are continuous, because $\varphi=\lambda^{-1}T_K\varphi$. By transitivity, $\operatorname{supp}\mu=\Omega$, so continuous functions that agree $\mu$-almost everywhere agree everywhere. We write $k_s(x,y)=\sum_{i\le s}\varphi_i(x)\varphi_i(y)$ for the reproducing kernel of $V_s$ in $L^2(\mu)$. Like every reproducing kernel of a finite-dimensional space, it does not depend on the choice of orthonormal basis \cite{Aronszajn1950}. We call $V_s$ \emph{$G$-invariant} if $\pi_gV_s=V_s$ for all $g$. This holds automatically when $\lambda_s>\lambda_{s+1}$, because $V_s$ is then a direct sum of full eigenspaces.

The main consequence of homogeneity used below is that the diagonals of the relevant kernels are constant. For Laplace eigenspaces on compact homogeneous Riemannian manifolds this is the addition formula of Gin\'e \cite{Gine1975}. Two further uses are more elementary. Theorem~\ref{thm:tightchar} uses that $\mu$ has full support, which follows from transitivity, and Theorem~\ref{thm:symmetry} uses the group action directly.

\begin{lemma}[Constant diagonal; cf.\ \cite{Gine1975}]\label{lem:diag}
In a homogeneous setting, suppose $V_s$ is $G$-invariant. Then for every $x\in\Omega$
\[
k_s(x,x)=s,\qquad K(x,x)=\sum_{i\ge1}\lambda_i,\qquad K^{(s)}(x,x)=\sum_{i\le s}\lambda_i,\qquad T^{(s)}(x,x)=\tau_s:=\sum_{i>s}\lambda_i .
\]
\end{lemma}
\begin{proof}
Let $W$ be a finite-dimensional $\pi_g$-invariant subspace of continuous functions with $L^2(\mu)$-orthonormal basis $w=(w_1,\dots,w_d)$ and reproducing kernel $k_W(x,y)=\sum_rw_r(x)w_r(y)$. Since $\mu$ is $G$-invariant, $\pi_gw$ is again an orthonormal basis of $W$, and basis independence of the reproducing kernel gives $k_W(g^{-1}x,g^{-1}y)=k_W(x,y)$ for $\mu\otimes\mu$-almost every $(x,y)$. Both sides are continuous and $\mu$ has full support, so the identity holds for all $x,y\in\Omega$. Apply this to $W=V_s$ for invariance of $k_s$. For $K^{(s)}$, note that $V_s$ is spanned by eigenfunctions, so $V_s=\bigoplus_\lambda(V_s\cap\mathcal E_\lambda)$, and each summand is $\pi_g$-invariant as an intersection of invariant subspaces. Choosing the basis of $V_s$ adapted to this decomposition, which is allowed by Remark~\ref{rem:basis}, gives $K^{(s)}=\sum_\lambda\lambda\,k_{V_s\cap\mathcal E_\lambda}$, a sum of invariant kernels. $K$ is invariant by hypothesis, and $T^{(s)}=K-K^{(s)}$. An invariant kernel $k$ has constant diagonal, since $k(gx,gx)=k(x,x)$ and $G$ acts transitively. The constants follow by integrating against $\mu$, \emph{i.e.} we have $\int k_s(x,x)\dd(x)=\sum_{i\le s}\|\varphi_i\|_{L^2}^2=s$ and $\int K^{(s)}(x,x)\dd(x)=\sum_{i\le s}\lambda_i$, while $\int K(x,x)\dd(x)=\sum_i\lambda_i$ by Mercer's theorem.
\end{proof}

\begin{remark}
For the Gaussian kernels on $\Sph^1$ and $\Sph^2$ we have $K(x,x)=1$, so $\tau_s=1-\sum_{i\le s}\lambda_i$, in agreement with Propositions~\ref{prop:circle} and~\ref{prop:sphere}.
\end{remark}

\subsection{Tight Sampling Sets}\label{sec:tight}

Let $\Phi\in\R^{m\times s}$ be the \emph{eigenfunction sampling matrix} with entries $\Phi_{ji}=\varphi_i(z_j)$ for $i\le s$, and let $\Lambda=\mathrm{diag}(\lambda_1,\dots,\lambda_s)$, so that $K^{(s)}_{ZZ}=\Phi\Lambda\Phi^\top$. We first take $m=s$, so that $\Phi$ is square.

In the orthonormal basis $\{\psi_i\}$ of $\Hcal$, the translate $K_{z_j}$ has leading coordinates $\sqrt{\lambda_i}\varphi_i(z_j)$ for $i\le s$, which form the $j$-th column of $D:=\Lambda^{1/2}\Phi^\top$. Its tail coordinates, taken in the completed basis of Section~\ref{sec:setting}, form the $j$-th column of an operator $D_\perp:\R^m\to V_s^\perp$ with $D_\perp^*D_\perp=T^{(s)}_{ZZ}$. Hence $K_{ZZ}=D^\top D+D_\perp^*D_\perp$. When $D$ is invertible, the tangents of the principal angles between $U_Z=\operatorname{range}[D;D_\perp]$ and $V_s$ are the singular values of $D_\perp D^{-1}$~\cite[Thm.~3.1]{ZhuKnyazev2013}. Invertibility of $D$ means that no principal angle equals $\pi/2$. Squaring these tangents gives the following reformulation of the closed form~\eqref{eq:closedform}.

\begin{proposition}[Tangent form of the alignment]\label{thm:factor}
Let $m=s$, suppose $\Phi$ is invertible, and set
 \[R_Z:=(D_\perp D^{-1})^*(D_\perp D^{-1})=\Lambda^{-1/2}\Phi^{-1}\,T^{(s)}_{ZZ}\,\Phi^{-\top}\Lambda^{-1/2}\succeq0 .\]
The eigenvalues of $R_Z$ are $\tan^2\theta_1,\dots,\tan^2\theta_s$, where $\theta_j$ are the principal angles between $U_Z$ and $V_s$. Consequently
\[
F(Z)=\sum_{j=1}^s\sin^2\theta_j=\tr\big((I+R_Z)^{-1}R_Z\big),\qquad \frac{\tr R_Z}{1+\|R_Z\|}\le F(Z)\le\tr R_Z .
\]
Moreover
\[
\tr R_Z=\sum_{j}\tan^2\theta_j=\sum_{i>s}\lambda_i\,\big\|\Lambda^{-1/2}c_i\big\|_2^2,\qquad c_i:=\Phi^{-1}\varphi_i(Z),
\]
where $c_i$ is the coefficient vector of the interpolant in $V_s$ of the tail eigenfunction $\varphi_i$ from its values at $Z$. We call $c_i$ the \emph{aliasing coefficients} of $\varphi_i$ onto $V_s$.
\end{proposition}
\begin{proof}
The spectral statement is \cite[Thm.~3.1]{ZhuKnyazev2013} applied to the bases above, and $\sin^2=\tan^2/(1+\tan^2)$ gives the formula for $F$, consistent with Section~\ref{sec:discrepancy}. Alternatively, $K_{ZZ}=D^\top(I+R_Z)D$ and \eqref{eq:closedform} give $F=\tr(K_{ZZ}^{-1}T^{(s)}_{ZZ})=\tr((I+R_Z)^{-1}R_Z)$ directly. The bounds follow from $t/(1+\|R_Z\|)\le t/(1+t)\le t$ for $t$ an eigenvalue of $R_Z$. The trace identity follows from expanding $T^{(s)}_{ZZ}=\sum_{i>s}\lambda_i\varphi_i(Z)\varphi_i(Z)^\top$.
\end{proof}

The spectral alignment is therefore squeezed between two saturating functions of the weighted aliasing energy $\tr R_Z=\|\tan\Theta\|_F^2$, and it is close to $\tr R_Z$ whenever all principal angles are small. The contribution that is specific to landmark selection is the last identity. Landmarks are good when the tail eigenfunctions, as seen through their samples at $Z$, look like low-energy elements of $V_s$. The next result turns this into a statement about which sampling matrices are best. The key is that $\|\Phi\|_F^2$ is fixed by Lemma~\ref{lem:diag}.

\begin{theorem}[Tight sampling sets minimize the bound]\label{thm:tight}
In a homogeneous setting with $V_s$ $G$-invariant and $m=s$, let
$\Phi\in\mathbb R^{s\times s}$, $\Phi_{ji}=\varphi_i(z_j)$, let
$\sigma_{\min}(\Phi)$ be its smallest singular value, and set
$\tau_s:=\sum_{i>s}\lambda_i$. Then
\[
F(Z)\ \le\ \tr R_Z\ \le\ \frac{s\,\tau_s}{\lambda_s\,\sigma_{\min}(\Phi)^2},\qquad\text{and}\qquad \sigma_{\min}(\Phi)^2\le s .
\]
Hence the bound is at least $\tau_s/\lambda_s$, with equality if and only if $\Phi^\top\Phi=sI$. We call such $Z$ a \emph{tight sampling set}.
\end{theorem}
\begin{proof}
If $\Phi$ is singular then $\sigma_{\min}=0$ and there is nothing to prove. Otherwise $D=\Lambda^{1/2}\Phi^\top$ is invertible, so
Proposition~\ref{thm:factor} applies and we have \[
  \operatorname{tr}R_Z
  =\sum_{i>s}\lambda_i\,\bigl\|\Lambda^{-1/2}\Phi^{-1}\varphi_i(Z)\bigr\|_2^2 .
\] The matrix $\Lambda^{-1/2}$ is diagonal with entries $\lambda_k^{-1/2}\le
\lambda_s^{-1/2}$ for $k\le s$, because the eigenvalues are ordered. Also
$\|\Phi^{-1}v\|_2\le\|v\|_2/\sigma_{\min}(\Phi)$. Therefore
\[
  \bigl\|\Lambda^{-1/2}\Phi^{-1}\varphi_i(Z)\bigr\|_2^2
  \le\frac{\|\varphi_i(Z)\|_2^2}{\lambda_s\,\sigma_{\min}(\Phi)^2},
\]
and summing,
\[
  \operatorname{tr}R_Z\le
  \frac{1}{\lambda_s\sigma_{\min}(\Phi)^2}\sum_{i>s}\lambda_i\|\varphi_i(Z)\|_2^2
  =\frac{\operatorname{tr}T^{(s)}_{ZZ}}{\lambda_s\sigma_{\min}(\Phi)^2}.
\] By Lemma~\ref{lem:diag}, $T^{(s)}(z_j,z_j)=\tau_s$ for every $j$, so
$\operatorname{tr}T^{(s)}_{ZZ}=s\,\tau_s$ and thus
\[
\tr R_Z\le\frac{s\,\tau_s}{\lambda_s\sigma_{\min}(\Phi)^2}.
\]
Let $\sigma_1,\dots,\sigma_s$ be the singular values of $\Phi$. Then
$\|\Phi\|_F^2=\sum_{k=1}^s\sigma_k^2$. Computing the same quantity entrywise
and applying Lemma~\ref{lem:diag} to $k_s$,
\[
  \|\Phi\|_F^2=\sum_{j=1}^s\sum_{i\le s}\varphi_i(z_j)^2
  =\sum_{j=1}^s k_s(z_j,z_j)=s\cdot s=s^2 .
\]
The smallest of $s$ nonnegative numbers is at most their mean, so
\[
  \sigma_{\min}(\Phi)^2\le\frac{1}{s}\sum_{k=1}^s\sigma_k^2=\frac{s^2}{s}=s .
\]
Equality holds if and only if $\sigma_k^2=s$ for all $k$, that is, if and
only if $\Phi^\top\Phi=sI$.
\end{proof}


The word ``tight'' is chosen deliberately. Below we show that tight sampling sets are tight spherical designs on $\Sph^2$ and unit-norm tight frames in $\R^s$. Their first intrinsic description comes from cubature. Write $V_s\cdot V_s=\Span\{fg:\ f,g\in V_s\}$. This space contains the constants, because $\sum_{i\le s}\varphi_i^2=s$.

\begin{theorem}[Tight sampling sets are minimal cubature rules]\label{thm:tightchar}
In a homogeneous setting with $V_s$ $G$-invariant, let $Z=\{z_1,\dots,z_s\}$ consist of $s$ distinct points. The following are equivalent.
\begin{enumerate}[label=(\roman*),leftmargin=*]
\item $\Phi^\top\Phi=sI$, so $Z$ is tight.
\item $k_s(z_j,z_l)=0$ for all $j\ne l$. Equivalently, the functions $k_s(\cdot,z_j)$ are pairwise orthogonal in $L^2(\mu)$.
\item The equal-weight rule $\frac{1}{s}\sum_jf(z_j)$ integrates every $f\in V_s\cdot V_s$ exactly.
\item Some rule $\sum_jw_jf(z_j)$ with positive weights integrates every $f\in V_s\cdot V_s$ exactly.
\item The Lagrange basis of $V_s$ at $Z$ exists and is $\ell_j=k_s(\cdot,z_j)/s$.
\end{enumerate}
\end{theorem}

The theorem is the homogeneous-space version of the classical theory of minimal (Gaussian) cubature \cite{Mysovskikh1981,CoolsMysovskikhSchmid2001}. A positive rule exact on $V\cdot V$ needs at least $\dim V$ nodes~\cite{Mysovskikh1981}, without any homogeneity assumption. A rule attaining this bound has nodes at which the reproducing kernel of $V$ vanishes pairwise, and weights equal to the Christoffel numbers $1/k_V(z_j,z_j)$. The only new ingredient is Lemma~\ref{lem:diag}, which makes the Christoffel function constant and hence forces equal weights $1/s$. Tight sampling sets are therefore exactly the minimal positive cubature rules for $V_s\cdot V_s$. On $\Sph^2$ with $s=(M+1)^2$, the space $V_s\cdot V_s$ consists of the polynomials of degree at most $2M$, so tight sampling sets are exactly the \emph{tight} spherical $2M$ -designs, namely those attaining the Delsarte--Goethals--Seidel (Fisher-type) bound $(M+1)^2$ \cite{DelsarteGoethalsSeidel1977}.

\begin{proof}
(i)$\Leftrightarrow$(ii): $\Phi$ is square, so $\Phi^\top\Phi=sI$ if and only if $\Phi\Phi^\top=sI$, and $(\Phi\Phi^\top)_{jl}=k_s(z_j,z_l)$ with diagonal $s$. (i)$\Leftrightarrow$(iii): $(\Phi^\top\Phi)_{ik}=\sum_j(\varphi_i\varphi_k)(z_j)$ while $\int\varphi_i\varphi_k\dd=\delta_{ik}$. (iii)$\Rightarrow$(iv) is immediate. For (iv)$\Rightarrow$(ii) we follow the classical argument \cite{CoolsMysovskikhSchmid2001}. Exactness on $f^2$ shows that $f\in V_s$ vanishing on $Z$ is zero, so the Lagrange basis $\ell_j$ exists (evaluation at $Z$ is bijective on $V_s$). Exactness on $\ell_j\ell_k$ gives $\int\ell_j\ell_k\dd=w_j\delta_{jk}$, so $\{\ell_j/\sqrt{w_j}\}$ is orthonormal and $k_s$ in this basis is $k_s(x,y) = \sum_j \ell_j(x)\ell_j(y)/w_j$. Hence $k_s(z_r,z_t)=\delta_{rt}/w_r$ and Lemma~\ref{lem:diag} then forces $w_r=1/s$. (ii)$\Leftrightarrow$(v): $k_s(\cdot,z_j)/s\in V_s$ takes the value $\delta_{jl}$ at $z_l$ if and only if (ii) holds.
\end{proof}

Theorem~\ref{thm:tightchar} reduces the existence of tight sets on each space to a concrete question. Are there $s$ points at which the single invariant function $k_s$ vanishes pairwise? Once $k_s$ is written down, condition (ii) settles the model cases:

\begin{corollary}[Circle]\label{cor:tight_circle}
On $\Sph^1$ with $s=2M+1$, the tight sampling sets are exactly the rotations of the $s$ equispaced points. Equivalently, the rotated trapezoidal rules are the only $s$-node positive rules exact for trigonometric polynomials of degree $2M$.
\end{corollary}

\begin{corollary}[Sphere]\label{cor:tight_sphere}
On $\Sph^2$ with $s=(M+1)^2$, the addition theorem gives $k_s(x,y)=p_M(x\cdot y)$ with $p_M=\sum_{n=0}^M(2n+1)P_n=(M+1)P_M^{(1,0)}$. So $Z$ is tight if and only if its pairwise inner products are roots of $p_M$. This is the Delsarte--Goethals--Seidel description of tight $2M$ -designs \cite{DelsarteGoethalsSeidel1977}. Tight $2M$-designs on $\Sph^2$ exist only for $M=1$ \cite{BannaiDamerell1979,BannaiDamerell1980}, where $p_1(t)=1+3t$ and the tight sets are exactly the regular tetrahedra. For $s\ge9$ no tight sampling set exists.
\end{corollary}

For $s\ge 9$ on the sphere, the optimizer must therefore trade tightness against the aliasing weights. Table~\ref{tab:tight} confirms numerically that the optimized sets come much closer to tightness than Fibonacci lattices do.

\subsection[Tight Frames and Fekete Points]{Tight Frames, Fekete Points, and the Case $m\ge s$}\label{sec:frames}

Tightness has two further classical readings, and the first extends to $m>s$. Write $u_j:=\Phi^\top e_j=(\varphi_1(z_j),\dots,\varphi_s(z_j))\in\R^s$ for the feature vector of the $j$-th landmark in $V_s$. By Lemma~\ref{lem:diag}, every $u_j$ has the same norm, $\|u_j\|^2=k_s(z_j,z_j)=s$. Their frame operator is $\sum_ju_ju_j^\top=\Phi^\top\Phi$, and their Gram matrix is $\Phi\Phi^\top=[k_s(z_j,z_l)]$. Following Benedetto and Fickus \cite{BenedettoFickus2003}, define the \emph{frame potential}
\[
\mathrm{FP}(Z):=\sum_{j,l=1}^m k_s(z_j,z_l)^2=\big\|\Phi^\top\Phi\big\|_F^2 .
\]
Recall that a determinant-maximizing set for an interpolation problem is called a \emph{Fekete point set} \cite{BosEtAl2010}, or an \emph{extremal fundamental system} on the sphere \cite{SloanWomersley2004}. Its natural least-squares analogue for $m>s$ is a \emph{D-optimal design}, which maximizes $\det(\Phi^\top\Phi)$ \cite{Bartel2025}. Both are well defined because a change of orthonormal basis of $V_s$ replaces $\Phi$ by $\Phi Q$ with $Q$ orthogonal.

\begin{proposition}[Tight frames, frame potential and Fekete points]\label{prop:frames}
In a homogeneous setting with $V_s$ $G$-invariant, let $m\ge s$. For every configuration $Z$ of $m$ points,
\[
\sigma_{\min}(\Phi)^2\le m,\qquad \det(\Phi^\top\Phi)\le m^s,\qquad \mathrm{FP}(Z)\ge m^2s .
\]
Equality in any one of the three holds if and only if $\Phi^\top\Phi=mI$. This in turn is equivalent to each of the following:
\begin{enumerate}[label=(\alph*),leftmargin=*]
\item $\{u_j/\sqrt s\}_{j=1}^m$ is a unit-norm tight frame for $\R^s$;
\item the equal-weight rule $\frac{1}{m}\sum_jf(z_j)$ integrates every $f\in V_s\cdot V_s$ exactly.
\end{enumerate}
For $m=s$, $\det(\Phi^\top\Phi)=\det(\Phi)^2$. So the tight sampling sets of Theorem~\ref{thm:tight} are exactly the configurations attaining the universal Fekete bound $|\det\Phi|\le s^{s/2}$, whenever such configurations exist. Their Lagrange functions satisfy $\sup_x|\ell_j(x)|=1$.
\end{proposition}
\begin{proof}
The eigenvalues $\sigma_k^2$ of $\Phi^\top\Phi$ are nonnegative with sum $\|\Phi\|_F^2=\sum_j k_s(z_j,z_j)=ms$. The minimum is at most the mean $m$. The arithmetic–geometric mean inequality gives $\prod_k\sigma_k^2\le m^s$. Cauchy--Schwarz gives $\mathrm{FP}(Z) = \sum_k\sigma_k^4\ge(ms)^2/s$, which is the Welch bound \cite{Welch1974}. Each inequality is an equality if and only if all $\sigma_k^2=m$. The frame operator of the normalized vectors $u_j/\sqrt s$ is $\Phi^\top\Phi/s$, so $\Phi^\top\Phi=mI$ says exactly that $\{u_j/\sqrt s\}$ is a tight frame with frame bound $m/s$, which is the bound forced on any unit-norm tight frame of $m$ vectors in $\R^s$ \cite{BenedettoFickus2003}. This is reformulation (a). Statement (b) is the computation (i)$\Leftrightarrow$(iii) of Theorem~\ref{thm:tightchar} with $m$ nodes. For the last claim, (v) of Theorem~\ref{thm:tightchar} and Cauchy--Schwarz give $|\ell_j(x)|=|k_s(x,z_j)|/s\le\sqrt{k_s(x,x)k_s(z_j,z_j)}/s=1$, with equality at $x=z_j$.
\end{proof}

For $m=s$, Proposition~\ref{prop:frames} shows that three different notions of good sampling coincide on tight sets: the conditioning criterion of Theorem~\ref{thm:tight} (maximal $\sigma_{\min}$), the Fekete criterion (maximal $|\det\Phi|$), and the frame-potential criterion (minimal $\sum_{j\neq l}k_s(z_j,z_l)^2$). When tight sets do not exist, as on $\Sph^2$ for $s\ge9$, the three criteria are genuinely different relaxations.

For $m>s$, Proposition~\ref{prop:frames} gives a clean extension of tightness. On $\Sph^1$ with $s=2M+1$, the $m$ equispaced points are a tight frame for every $m\ge s$, because the trapezoidal rule with $m$ nodes is exact up to degree $m-1\ge 2M$. On $\Sph^2$ with $s=(M+1)^2$, the condition $\Phi^\top\Phi=mI$ says exactly that $Z$ is a spherical $2M$ -design with $m$ points. For $s=4$, examples are the octahedron and the cube (both $3$-designs) and the icosahedron and dodecahedron (both $5$-designs), and the last two also qualify for $s=9$. The equivalence with frames is the familiar correspondence between designs and tight frames \cite{BenedettoFickus2003,Bartel2025}.

What does \emph{not} extend is Proposition~\ref{thm:factor}, because for $m>s$ the principal angles are no longer controlled by $\Phi$ alone. Tightness still describes the best-conditioned sampling of $V_s$, since it maximizes $\sigma_{\min}(\Phi)^2$ and $\det(\Phi^\top\Phi)$, but we have no analogue of Theorem~\ref{thm:tight} that links this conditioning to $F$ when $m>s$. The cube is a tight frame for $s=4$, yet for $s=4$ and $m=8$ it is a saddle point of $F$ (Table~\ref{tab:poly}), and the optimizer prefers the square antiprism.

\subsection{A Symmetry Argument}\label{sec:symmetry}

Tightness explains why even spacing is favorable, but it does not by itself show that symmetric configurations are stationary points of the objective. That follows from a classical principle of Michel \cite{Michel1971,MichelZhilinskii2001}. For a smooth action of a compact (in particular finite) group on a finite-dimensional manifold, the gradient of \emph{every} invariant function vanishes on orbits that are isolated in their stratum. For the action of $\Gamma\times S_m$ on landmark configurations, this hypothesis reduces to a condition on the stabilizers of the individual landmarks.

\begin{theorem}[Symmetric orbits are critical points]\label{thm:symmetry}
Let $\Omega$ be a Riemannian manifold, let $K$ be $C^1$, and let $\Gamma$ be a finite group of isometries preserving $\mu$ and $K$. Suppose $V_s$ is $\Gamma$-invariant (for example, a union of full eigenspaces). Let $Z=\{z_1,\dots,z_m\}$ be a $\Gamma$-orbit such that, for every $j$, the stabilizer $\Gamma_{z_j}$ fixes no nonzero tangent vector in $T_{z_j}\Omega$. Then $Z$ is a critical point of $F$, and more generally of every $C^1$ function on $\Omega^m$ that is invariant under permutations of the landmarks and under the diagonal action of $\Gamma$.
\end{theorem}
\begin{proof}
$F$ is invariant under both actions because $K_{ZZ}$ and $K^{(s)}_{ZZ}$ are. The subgroup $H=\{(g,\sigma_g):g\in\Gamma\}\subset\Gamma\times S_m$, where $\sigma_g$ is the permutation induced by $g$ on $Z$, fixes the ordered configuration $Z$. A tangent vector $v=(v_j)$ at $Z$ fixed by $H$ satisfies $dg\,v_j=v_{\sigma_g(j)}$. For $g\in\Gamma_{z_j}$ this says that $v_j$ is fixed by $\Gamma_{z_j}$, so $v_j=0$. Thus $Z$ is an isolated fixed point of $H$, i.e.\ its orbit is isolated in its stratum, and Michel's theorem \cite{Michel1971} applies. Directly, the gradient of an $H$-invariant function at an $H$-fixed point is an $H$-fixed tangent vector, hence zero.
\end{proof}

On the sphere, configurations whose stabilizers fix no tangent vector are what Cohn, Elkies, Kumar and Sch\"urmann call \emph{group-balanced} \cite{CohnElkiesKumarSchurmann2010}. Leech classified the configurations on $\Sph^2$ that are in equilibrium under every pairwise force law \cite{Leech1957}. Theorem~\ref{thm:symmetry} is the same principle applied to an objective that is not a pairwise energy. On the circle, the $m$ equispaced points form an orbit of the dihedral group $D_m$, and the reflection through each point acts as $-1$ on its tangent line. So equispaced landmarks are critical for every $m$, $s$. On the sphere, the vertices of the five Platonic solids have stabilizers containing a rotation of order $3$, $4$, or $5$ about the vertex axis, which fixes no tangent vector. All five solids are therefore critical points, and numerically $\|\nabla F\|\le3\times10^{-14}$ at each of them (Table~\ref{tab:poly}).

Because Theorem~\ref{thm:symmetry} applies to every invariant objective, it carries no information specific to $F$, and in particular none about minimality. The Hessian at the cube and at the dodecahedron has negative directions (Table~\ref{tab:poly}). This mirrors the fact that the cube is not an optimal eight-point configuration for packing or energy problems either. For $s=4$ and $m=8$ the optimizer instead finds a square antiprism. It has the combinatorial type of the best packing of $8$ points~\cite{SchutteVanDerWaerden1951}, but it is not that packing itself, because its edges have two different lengths (Section~\ref{sec:sphere}).

\subsection{Equispaced Landmarks on the Circle}\label{app:alias}

On the circle we can compute the leakage of the equispaced configuration exactly, for any number of landmarks $m\ge s$, and even the individual principal angles. The Gram matrices of equispaced points are circulant, so the discrete Fourier transform diagonalizes them \cite{Davis1979}. Sampling then aliases frequency $n$ onto the residue $n\bmod m$. Both facts are standard. The contribution here is the resulting closed form for $F$. Here and in Appendix~\ref{app:circlesphere} we index the circle eigenvalues by frequency, writing $\beta_k=e^{-\alpha}I_k(\alpha)$ for the eigenvalue of $e^{\pm ik\theta}$, so that the sorted eigenvalues are $\lambda_1=\beta_0$ and $\lambda_{2k}=\lambda_{2k+1}=\beta_k$ for $k\ge1$.

\begin{theorem}[Aliasing formula]\label{thm:alias}
Let $s=2M+1$, $m\ge s$, and let $Z$ be $m$ equispaced points on $\Sph^1$. The principal angles between $U_Z$ and $V_s$ are indexed by the frequencies $|k|\le M$, with
\[
\sin^2\theta_k=\frac{\sum_{j\ne0}\beta_{|k+jm|}}{\sum_{j\in\mathbb Z}\beta_{|k+jm|}},\qquad\text{so}\qquad
F(Z)=\sum_{k=-M}^{M}\frac{\sum_{j\ne0}\beta_{|k+jm|}}{\sum_{j\in\mathbb Z}\beta_{|k+jm|}} .
\] 
The condition number of the Gram matrix is
 \[
\kappa(K_{ZZ})=\frac{\max_r\sum_j\beta_{|r+jm|}}{\min_r\sum_j\beta_{|r+jm|}} .
\]
For $m=s$, the matrix $R_Z$ of Proposition~\ref{thm:factor} is diagonal in the Fourier basis with entries $\tan^2\theta_k=\sum_{j\ne0}\beta_{|k+jm|}/\beta_{|k|}$.
\end{theorem}

\begin{figure}[h!]
\centering
\includegraphics[width=\linewidth]{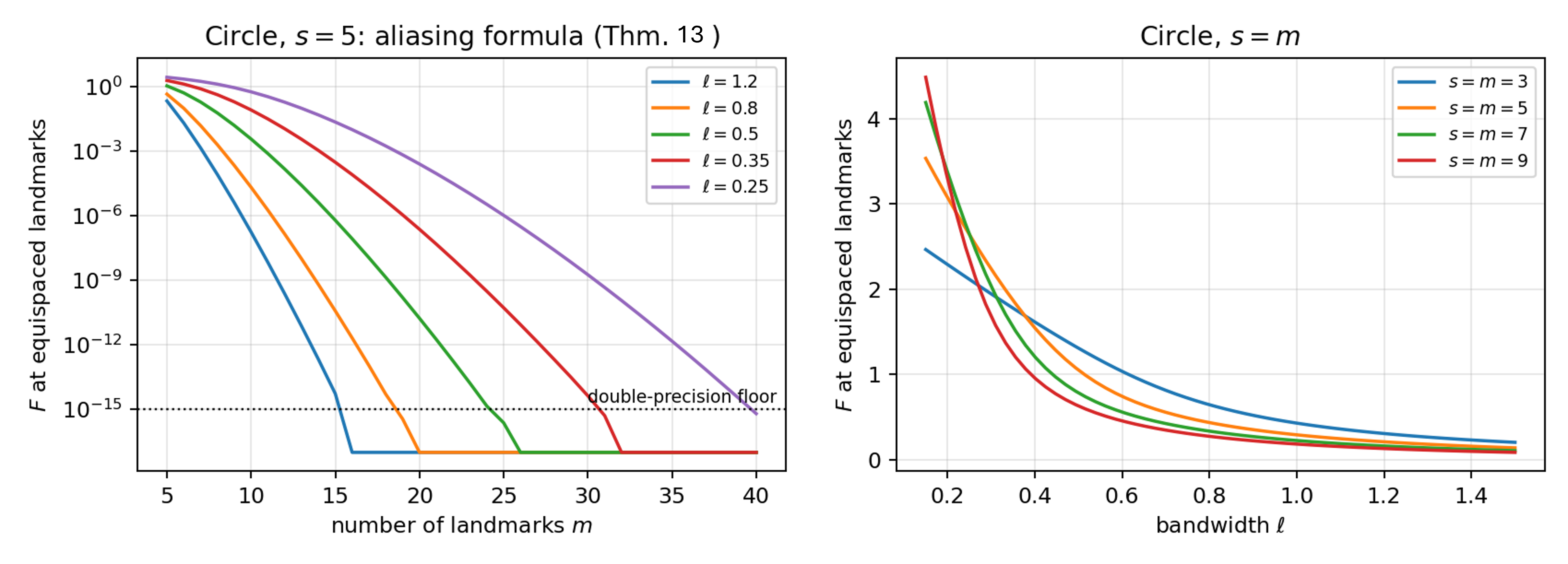}
\caption{Theorem~\ref{thm:alias} on the circle. Left: $F$ at $m$ equispaced landmarks for $s=5$ and several bandwidths. The decay is super-exponential, and the dotted line marks the double-precision resolution of $F=s-\tr(\cdot)$, with values below $10^{-16}$ clipped. Right: $F$ at equispaced landmarks with $m=s$, as a function of $\ell$.}
\label{fig:alias}
\end{figure}

See Appendix~\ref{app:circlesphere} for the proof. For large $m$ the formula shows that the dominant alias of the highest retained frequency $M$ is $m-M$, and the two frequencies $\pm M$ contribute equally, so $F\approx2\beta_{m-M}/\beta_M$. Since $\beta_k\sim e^{-\alpha}(\alpha/2)^k/k!$, the leakage decays super-exponentially in $m$ (Figure~\ref{fig:alias}).

\section{The Algorithm}\label{sec:algorithm}

The algorithm treats the landmark coordinates $Z=(z_1,\dots,z_m)$ as free
parameters and minimizes a differentiable, Cholesky-based form of the leakage
objective by gradient descent. It begins with a one-time precomputation of the top-$s$ eigenpairs of the data kernel matrix $K_N\approx U\Sigma U^\top$, with
$U\in\R^{N\times s}$ and $\Sigma=\mathrm{diag}(\nu_1,\dots,\nu_s)$. The landmarks are
then updated iteratively from an initial configuration until the gradient or change in $F$ reaches the
resolution of the objective.

To allow the landmarks to leave the data points, we extend the eigenvectors to any $z\in\R^d$ via the Nystr\"om extension~\cite{bengio2004}: \[
\psi_i(z)=\frac{1}{\sqrt{\nu_i}}\sum_{j=1}^N k(z,x_j)\,U_{ji},
\qquad\text{that is,}\qquad
\psi(z)^\top=k(z,X)\,U\Sigma^{-1/2}\in\R^{1\times s}.
\] 
Each $\psi_i=\sum_jc_jk(\cdot,x_j)$ has coefficient vector $c=u_i/\sqrt{\nu_i}$, so its squared RKHS norm is $c^\top K_Nc=u_i^\top K_Nu_i/\nu_i=1$, and distinct $\psi_i$ are orthogonal. At a data point $z=x_j$ the formula returns $\sqrt{\nu_i}\,U_{ji}$. Consequently $V_s=\Span\{\psi_i\}$ is an exact $s$-dimensional subspace of $\Hcal$, and all the projector algebra of the earlier sections applies unchanged, except that the points $z_j$ may now lie anywhere and not only at data points. These $\psi_i$ are exactly the functions $\sqrt{\lambda_i}\varphi_i$ of Section~\ref{sec:setting} for the empirical measure $\mu_N$, with $\lambda_i=\nu_i/N$. Because $\mu_N$ has finite support, the translates $K_{z_j}$ of landmarks off the data have components in $\Hcal_0$, and it is the completed basis of Section~\ref{sec:setting} that keeps the projector algebra, the mass conservation \eqref{eq:massconservation}, and the bounds of Section~\ref{sec:sandwich} valid in this case. 

For $m$ landmarks, define the matrices $H = U\Sigma^{-1/2}\in \R^{N\times s}$ and $\Psi(Z)=k(Z,X)\,H\in\R^{m\times s}$, whose $(j,i)$ entry is $\psi_i(z_j)$. For the purpose of ensuring a stable Cholesky factorization, also define $A(Z)=k(Z,Z)+\varepsilon I\in\R^{m\times m}$, where $\varepsilon$ is a small jitter that keeps $A$ numerically positive definite. By \eqref{eq:overlap}, $\sum_{i\le s}a_i=\tr(A^{-1}\Psi\Psi^\top)$. Taking the Cholesky factor $A=LL^\top$ and setting $Y=L^{-1}\Psi$ gives
\begin{equation}\label{eq:Fchol}
\boxed{F(Z)=s-\|L^{-1}\Psi(Z)\|_F^2}
\end{equation}
The Frobenius-norm form avoids forming $A^{-1}$. It needs one Cholesky factorization and one triangular solve, both of which autodiff frameworks differentiate stably. Since $\Psi\Psi^\top=k(Z,X)\,U\Sigma^{-1}U^\top\,k(X,Z)$, the objective can equivalently be read as a comparison between the landmark Gram matrix and the Nystr\"om-extended rank-$s$ kernel evaluated at $Z$.

\begin{algorithm}[h!]
\caption{Spectral alignment landmark optimization}\label{alg:landmarks}
\begin{algorithmic}[1]
\Require Data $X=\{x_i\}_{i=1}^N$, kernel $k$, target rank $s$, number of landmarks $m\ge s$, jitter $\varepsilon$, gradient tolerance $\mathrm{tol}_g$, relative decrease tolerance $\mathrm{tol}_F$, patience $p$
\Ensure Landmarks $Z=(z_1,\dots,z_m)$
\State Form $K_N$ and compute its top-$s$ eigenpairs $(U,\Sigma)$ by Lanczos or randomized SVD
\State $H\gets U\Sigma^{-1/2}$
\State Initialize $Z^{(0)}$ by greedy selection, random sampling, or $k$-means++
\State $c\gets0$
\For{$t=0,1,2,\dots$}
\State $\Psi\gets k(Z^{(t)},X)\,H$ \Comment{Nystr\"om-extended eigenfunctions}
\State $A\gets k(Z^{(t)},Z^{(t)})+\varepsilon I$
\State $L\gets\mathrm{chol}(A)$, $\quad Y\gets L^{-1}\Psi$
\State $F^{(t)}\gets s-\|Y\|_F^2$
\State $g\gets\nabla_ZF^{(t)}$ by reverse-mode autodiff
\If{$\|g\|_\infty<\mathrm{tol}_g$}
\State \Return $Z^{(t)}$
\EndIf
\State $Z^{(t+1)}\gets$ Adam or L-BFGS step on $Z^{(t)}$ using $g$
\State Retract $Z^{(t+1)}$ onto the manifold if constrained \Comment{e.g.\ $z_j\gets z_j/\|z_j\|$ on $\Sph^2$}
\EndFor
\end{algorithmic}
\end{algorithm}

On $\Sph^1$ we parametrize landmarks by angles. On $\Sph^2$ we use unconstrained vectors $v_j\in\R^3$ with $z_j=v_j/\|v_j\|$ and renormalize after every optimizer step, which avoids the pole singularities of latitude--longitude coordinates and makes the Euclidean gradient tangent. The objective is smooth and deterministic, so L-BFGS usually works well, especially in the badly conditioned tail.

\subsection{The gradient written out}

Autodiff computes the gradient automatically, but the explicit formula is useful for checking an implementation and for seeing how landmarks interact. We derive it for the data-driven objective. Throughout, $k:=k(Z,X)\in\R^{m\times N}$, $A=k(Z,Z)+\varepsilon I$, and $G=U\Sigma^{-1}U^\top$. Since $\Psi=kH=kU\Sigma^{-1/2}$, we have
\[
J:=kGk^\top=\Psi\Psi^\top,\qquad F=s-L_{\mathrm{tr}},\qquad L_{\mathrm{tr}}:=\tr(A^{-1}J).
\]
The landmarks enter $L_{\mathrm{tr}}$ through $A$, and hence $A^{-1}$, and through $k$, and hence $J$. The matrix $G$ depends only on the data and is constant.

By the product rule, $dL_{\mathrm{tr}}=\tr(dA^{-1}\,J)+\tr(A^{-1}\,dJ)$. Differentiating $A^{-1}A=I$ gives $dA^{-1}=-A^{-1}\,dA\,A^{-1}$, and since $G$ is constant, $dJ=dk\,G\,k^\top+k\,G\,dk^\top$. Hence
\[
dL_{\mathrm{tr}}=-\tr(A^{-1}\,dA\,A^{-1}\,J)+\tr(A^{-1}\,dk\,G\,k^\top)+\tr(A^{-1}\,k\,G\,dk^\top).
\]

Transposing the last term inside the trace and using the symmetry of $A^{-1}$ and $G$ gives $\tr(A^{-1}\,k\,G\,dk^\top)=\tr(dk\,G\,k^\top A^{-1})=\tr(A^{-1}\,dk\,G\,k^\top)$, so the two $dk$ terms coincide and together equal $2\tr(G\,k^\top A^{-1}\,dk)=\langle2A^{-1}kG,dk\rangle$. Here we write each term as a Frobenius inner product $\langle\Delta,dx\rangle=\tr(\Delta^\top dx)$ for the matrices $\Delta$ defined below. By cyclicity, the first term is $-\tr(A^{-1}JA^{-1}\,dA)=\langle-A^{-1}JA^{-1},dA\rangle$, because $A^{-1}JA^{-1}$ is symmetric. Therefore
\[
dL_{\mathrm{tr}}=\langle\Delta_k,dk\rangle+\langle\Delta_A,dA\rangle,\qquad \Delta_k=2\,A^{-1}\,k\,G\in\R^{m\times N},\qquad \Delta_A=-\,A^{-1}\,J\,A^{-1}\in\R^{m\times m}.
\]

Write $\nabla_1k$ and $\nabla_2k$ for the gradients of the kernel in its first and second arguments. The entry $k_{jl}=k(z_j,x_l)$ depends only on $z_j$, with $dk_{jl}=\nabla_1k(z_j,x_l)\cdot dz_j$. The entry $A_{jl}=k(z_j,z_l)+\varepsilon\delta_{jl}$ depends on $z_j$ through row $j$ and through column $j$. Collecting the terms that contain $dz_j$ gives
\[\sum_{l=1}^m(\Delta_A)_{jl}\,\nabla_1k(z_j,z_l)+\sum_{l=1}^m(\Delta_A)_{lj}\,\nabla_2k(z_l,z_j)=2\sum_{l=1}^m(\Delta_A)_{jl}\,\nabla_1k(z_j,z_l),
\] 
where the equality uses $\nabla_2k(z_l,z_j)=\nabla_1k(z_j,z_l)$, which follows from the symmetry of $k$, and the symmetry of $\Delta_A$. For the diagonal entry $l=j$, observe that $\frac{d}{dz}k(z,z)=2\nabla_1k(z,z)$, and this term vanishes for any stationary kernel. Combining with the $dk$ term and with $F=s-L_{\mathrm{tr}}$ gives
 \begin{equation}\label{eq:gradF}
\nabla_{z_j}F=-\Big[\sum_{l=1}^N(\Delta_k)_{jl}\,\nabla_1k(z_j,x_l)\;+\;2\sum_{l=1}^m(\Delta_A)_{jl}\,\nabla_1k(z_j,z_l)\Big].
\end{equation}

The first sum in \eqref{eq:gradF} couples each landmark to the data, weighted by how the extended eigenfunctions $\psi_i$ load on each data point after preconditioning by $A^{-1}$. The second sum couples the landmarks to one another. 
\section{Numerical Conditioning}\label{sec:conditioning}

The Gaussian kernel $k(x,y)=\exp(-\|x-y\|^2/2\ell^2)$ has a single length
scale $\ell$, the bandwidth, which sets how far a translate $k(\cdot,z_j)$
reaches. For large $\ell$ the translates are broad and nearly parallel, and for
small $\ell$ they are narrow and nearly orthogonal. This governs the
conditioning of the landmark Gram matrix $K_{ZZ}$, which the Cholesky-based
objective inverts at every iteration. Broad translates make the columns of
$K_{ZZ}$ nearly dependent and $\kappa(K_{ZZ})$ large, while narrow translates
make $K_{ZZ}$ close to the identity.

Consider the landmark spacing $h$, which is
$2\pi/m$ on $\Sph^1$ and of order $\sqrt{4\pi/m}$ on $\Sph^2$~\cite{SaffKuijlaars1997}. Adding landmarks decreases $h$, so at fixed $\ell$ more translates overlap and
the Gram matrix becomes worse conditioned. Keeping the Cholesky-based objective numerically well posed therefore requires a smaller bandwidth as $m$ grows. This section first quantifies this conditioning limit, which is an upper limit on $\ell$ for a given $m$, and then identifies an opposite lower limit, below which the objective depends on the landmark positions only through terms too small to resolve in floating point.

\begin{figure}[h!]
\centering
\includegraphics[width=\linewidth]{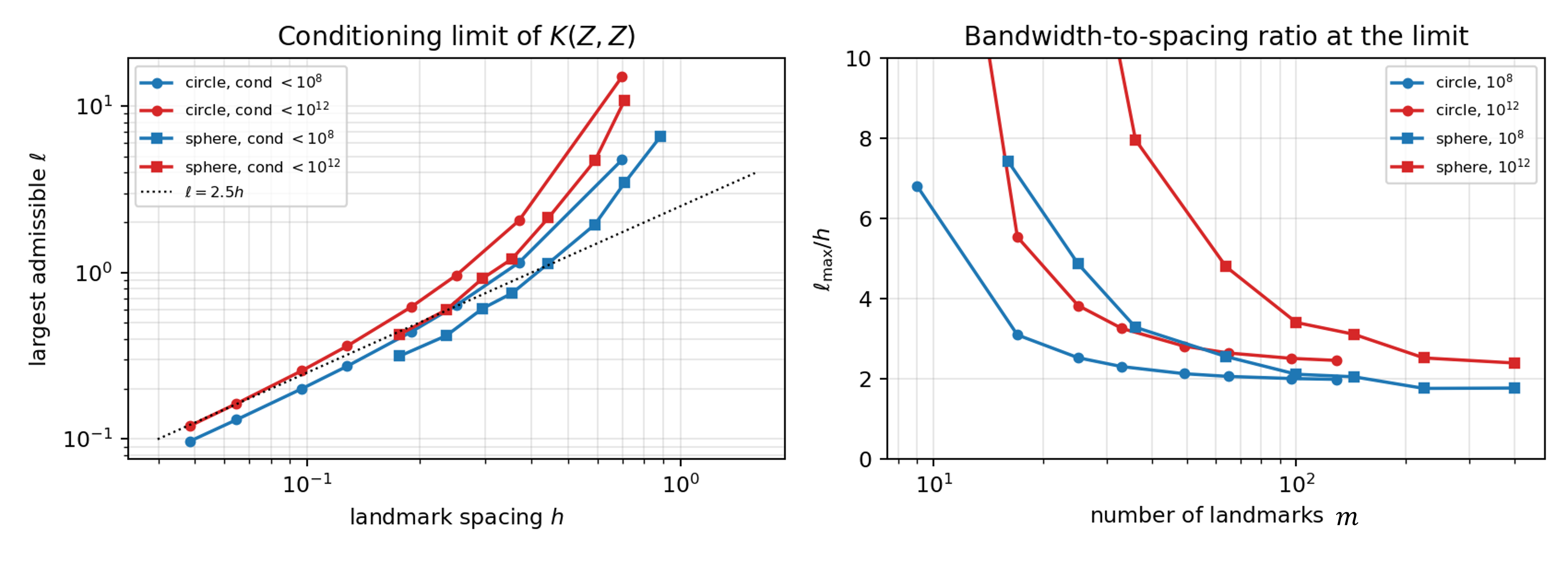}
\caption{Conditioning limit of $K(Z,Z)$. Left: largest admissible bandwidth against landmark spacing, with the reference line $\ell=2.5h$. Right: the ratio $\ell_{\max}/h$, which settles to a constant on $\Sph^1$ and decreases slowly toward a somewhat smaller value on $\Sph^2$.}
\label{fig:cond}
\end{figure}

On $\Sph^1$ and $\Sph^2$ the Gaussian
eigenvalues decay super-exponentially~\cite{NISTHandbook2010, SuperExpDecaySphere2020}, so in floating point the Gram matrix has a finite numerical rank. Landmark counts much larger than this rank therefore force $K_{ZZ}$ to have
eigenvalues below the working precision, and so $\kappa(K_{ZZ})$ becomes
prohibitively large. The rank grows as the bandwidth shrinks, since a narrower kernel has a slower-decaying spectrum and more eigenvalues remain
above the working precision. Solving $\kappa(K_{ZZ})=10^{12}$ for $\ell$ (Figure~\ref{fig:cond}) gives the largest admissible bandwidths in Table~\ref{tab:cond}. Beyond this limit two things can go wrong. The Cholesky factorization either fails or is dominated by the jitter $\varepsilon$, which biases $F$ because $\varepsilon$ then exceeds the smallest eigenvalues of $K_{ZZ}$. In addition, errors in $\nabla\tr(A^{-1}K^{(s)}_{ZZ})$ are amplified by $\kappa(A)$, so the descent direction becomes noise. This is Schaback's uncertainty principle \cite{Schaback1995} at work, since flat kernels approximate smooth functions well but produce ill-conditioned systems. If a larger $\ell$ is genuinely needed, a stable basis should be constructed. With the analytic expansion of Appendix~\ref{sec:homog}, one can factor the decaying $\beta_n$ out of $K_{ZZ}$ analytically, in the spirit of RBF-QR \cite{FornbergPiret2008,FornbergLarssonFlyer2011}.

\begin{table}[h]
\centering
\begin{tabular}{@{}rcc@{\qquad}rcc@{}}
\toprule
\multicolumn{3}{c}{$\Sph^1$, equispaced, $h=2\pi/m$} & \multicolumn{3}{c}{$\Sph^2$, Fibonacci, $h=\sqrt{4\pi/m}$}\\
$m$ & $\ell_{\max}$ & $\ell_{\max}/h$ & $m$ & $\ell_{\max}$ & $\ell_{\max}/h$\\
\midrule
17 & 2.05 & 5.55 & 36 & 4.71 & 7.97\\
33 & 0.62 & 3.26 & 64 & 2.13 & 4.80\\
49 & 0.36 & 2.82 & 100 & 1.21 & 3.41\\
65 & 0.26 & 2.65 & 144 & 0.92 & 3.12\\
97 & 0.16 & 2.52 & 225 & 0.60 & 2.53\\
129 & 0.12 & 2.47 & 400 & 0.43 & 2.40\\
\bottomrule
\end{tabular}
\caption{Largest Gaussian bandwidth keeping $\kappa(K_{ZZ})<10^{12}$. On $\Sph^1$ the ratio $\ell_{\max}/h$ converges to about $2.4$ (it is $2.41$ at $m=257$ and $2.40$ at $m=513$), so $\ell_{\max}\propto m^{-1}$. On $\Sph^2$ the ratio is still decreasing at $m=400$. Recomputing it for Fibonacci sets with $m=900$, $1600$ and $2500$ gives values between $1.9$ and $2.2$, so $\ell_{\max}\propto m^{-1/2}$ with a constant somewhat below the circle value. With the threshold $10^{8}$ the ratios are about $2.0$ on $\Sph^1$ and between $1.5$ and $1.8$ on $\Sph^2$ over the same range of $m$.}
\label{tab:cond}
\end{table}

The bandwidth is also bounded from below. If $\ell\ll q_Z$, where $q_Z$ is half the minimum pairwise distance of the
landmarks, then the kernels are too narrow to see each other and
$K_{ZZ}\approx I$. In that case write $K_{ZZ}=I+O_{ZZ}$, where $O_{ZZ}$ has zero diagonal and its off-diagonal entries are $e^{-\|z_j-z_l\|^2/2\ell^2}$. Whenever $\|O_{ZZ}\|<1$, the Neumann series $(I+O_{ZZ})^{-1}=I-O_{ZZ}+O_{ZZ}^2-\cdots$ converges, and substituting it into~\eqref{eq:F} gives
\[F(Z)=\underbrace{s-\tr K^{(s)}_{ZZ}}_{\text{leading term}}\;+\;\underbrace{\tr\big(O_{ZZ}\,K^{(s)}_{ZZ}\big)-\tr\big(O_{ZZ}^2K^{(s)}_{ZZ}\big)+\cdots}_{\text{corrections}}.\]
By Lemma~\ref{lem:diag}, the leading term is $s-\tr K^{(s)}_{ZZ}=s-m(1-\tau_s)$, which does not depend on the positions. The first position-dependent term is $\tr(O_{ZZ}K^{(s)}_{ZZ})$. Since every pairwise distance is at least $2q_Z$, each entry of $O_{ZZ}$ is at most $e^{-2q_Z^2/\ell^2}$, and the entries between nearest neighbors attain this order. For equispaced points ($q_Z=h/2$) at $\ell=0.15h$ this factor is $e^{-1/(2\cdot0.15^2)}\approx2\times10^{-10}$. The landmark positions therefore still determine $F$, but only through a correction about ten orders of magnitude below the leading term, which the optimizer cannot resolve with absolute tolerances in double precision. The lower limit is thus a limit of numerical resolution, not a loss of information in the objective.

\begin{figure}[h!]
\centering
\includegraphics[width=\textwidth]{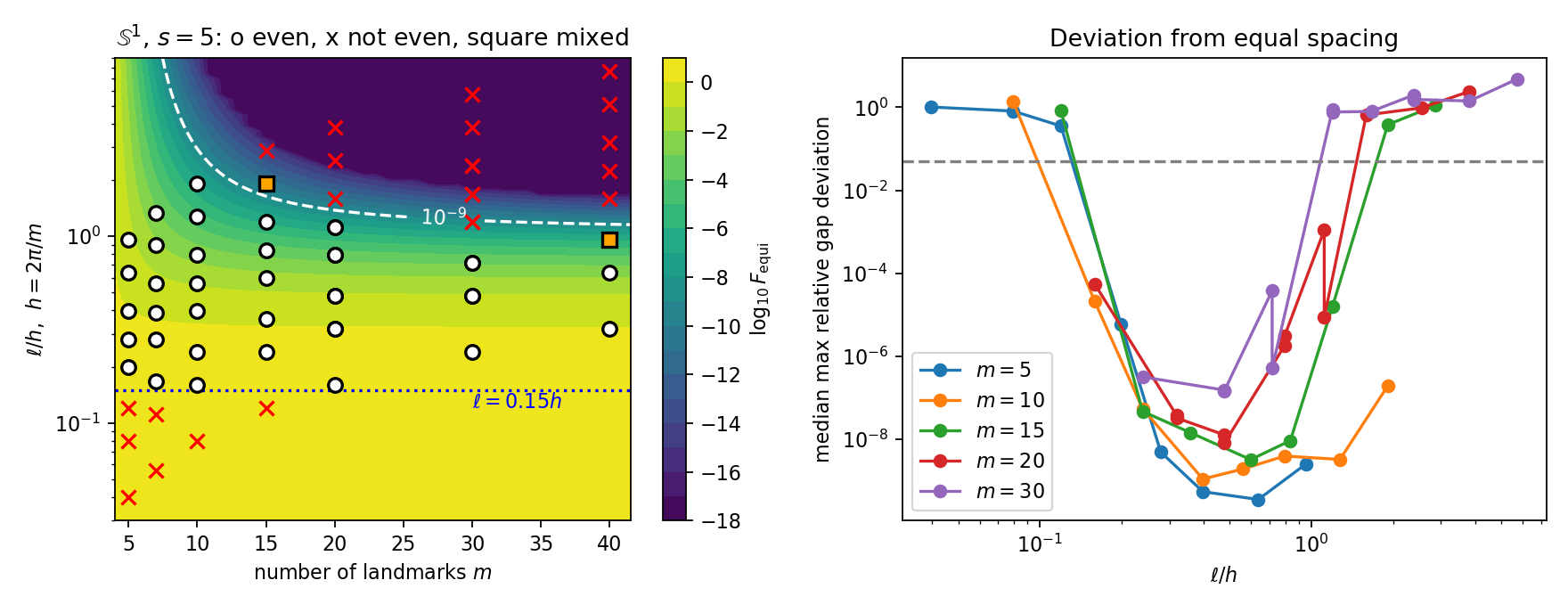}
\caption{Trade-off between $m$ and $\ell$ on $\Sph^1$, $s=5$. Left: markers show the outcome of the multistart runs in the $(m,\ell/h)$ plane (circle: all runs evenly spaced, square: some, cross: none). The background is $\log_{10}F_{\rm equi}(m,\ell)$ from Theorem~\ref{thm:alias}. The dashed white line is $F_{\rm equi}=10^{-9}$ and the dotted line is $\ell=0.15h$. Right: median of the largest relative gap deviation versus $\ell/h$. The dashed line is the $5\%$ evenness threshold.}
\label{fig:tradeoff_circle}
\end{figure}

To locate the two limits empirically, we swept $m$ and $\ell$ for $s=5$ on $\Sph^1$ with the analytic objective and uniform random initialization. We recorded whether the result was evenly spaced, meaning that the maximum relative deviation of consecutive gaps from $2\pi/m$ is below $5\%$. Figure~\ref{fig:tradeoff_circle} shows the outcome.

The results support two conclusions. First, the lower limit is sharp and independent of $m$. Every failure below the window occurs at $\ell/h\le0.15$, which is consistent with the flat-landscape analysis above. 
 Second, the upper limit decreases as $m$ grows, so the window shrinks. The admissible upper value of $\ell/h$ therefore decreases from $\ge1.9$
at $m=10$ to around $0.98$ at $m=30$, while the lower limit stays
near $0.15$. Runs above the window fail once $F_{\rm equi}$ (Theorem~\ref{thm:alias}) drops below about
$10^{-9}$, even when the Gram matrix is still well conditioned. At $m=30$ and
$\ell=0.25$, for example, $\kappa\approx440$ but $F_{\rm equi}=1.8\times10^{-9}$,
and the optimizer stops at uneven configurations with $F\approx7\times10^{-8}$.
Once the minimum is that small, the objective and its gradient are small as
well, so many uneven configurations satisfy the optimizer's absolute
stopping tolerances even though their $F$ is still above the minimum. At still larger $\ell$ (with
$\kappa\gtrsim10^{8}$), ill-conditioning takes over and the computed $F$
even becomes negative, with a median of $-8.5\times10^{-3}$ at $m=20$ and
$\ell=1.2$. Inside the window, the recovered configurations are equispaced to
high accuracy, with median gap deviations between $10^{-10}$ and $10^{-5}$,
and the attained $F$ agrees with $F_{\rm equi}$ to at least six digits.

\begin{remark}[Resolving small values of $F$]\label{rem:smallF}
Both failure modes at the edges of the window are failures to resolve a small objective, and they are aggravated by the way $F$ is computed. In the Cholesky form $F=s-\|Y\|_F^2$ the objective is the difference of two numbers of size $s$, so its absolute rounding error is at least of order $s\,\epsilon_{\rm mach}$, grows with $\kappa(K_{ZZ})$, and becomes large relative to $F$ as $F\to0$. This is also how negative values of $F$ arise. For $m=s$, where $K^{(s)}_{ZZ}$ is invertible, the identity $F=\tr\big(K_{ZZ}^{-1}T^{(s)}_{ZZ}\big)$ of Section~\ref{sec:discrepancy} (equation~\eqref{eq:F}) expresses $F$ through the tail kernel directly, and it avoids the cancellation. Independently of how $F$ is evaluated, minimizing $\log F$ or using stopping tolerances relative to the current value of $F$ removes the dependence of the stopping rule on the size of the minimum. We have not rerun the sweeps of this section with these modifications.
\end{remark}

\begin{figure}[h!]
\centering
\includegraphics[width=0.7\textwidth]{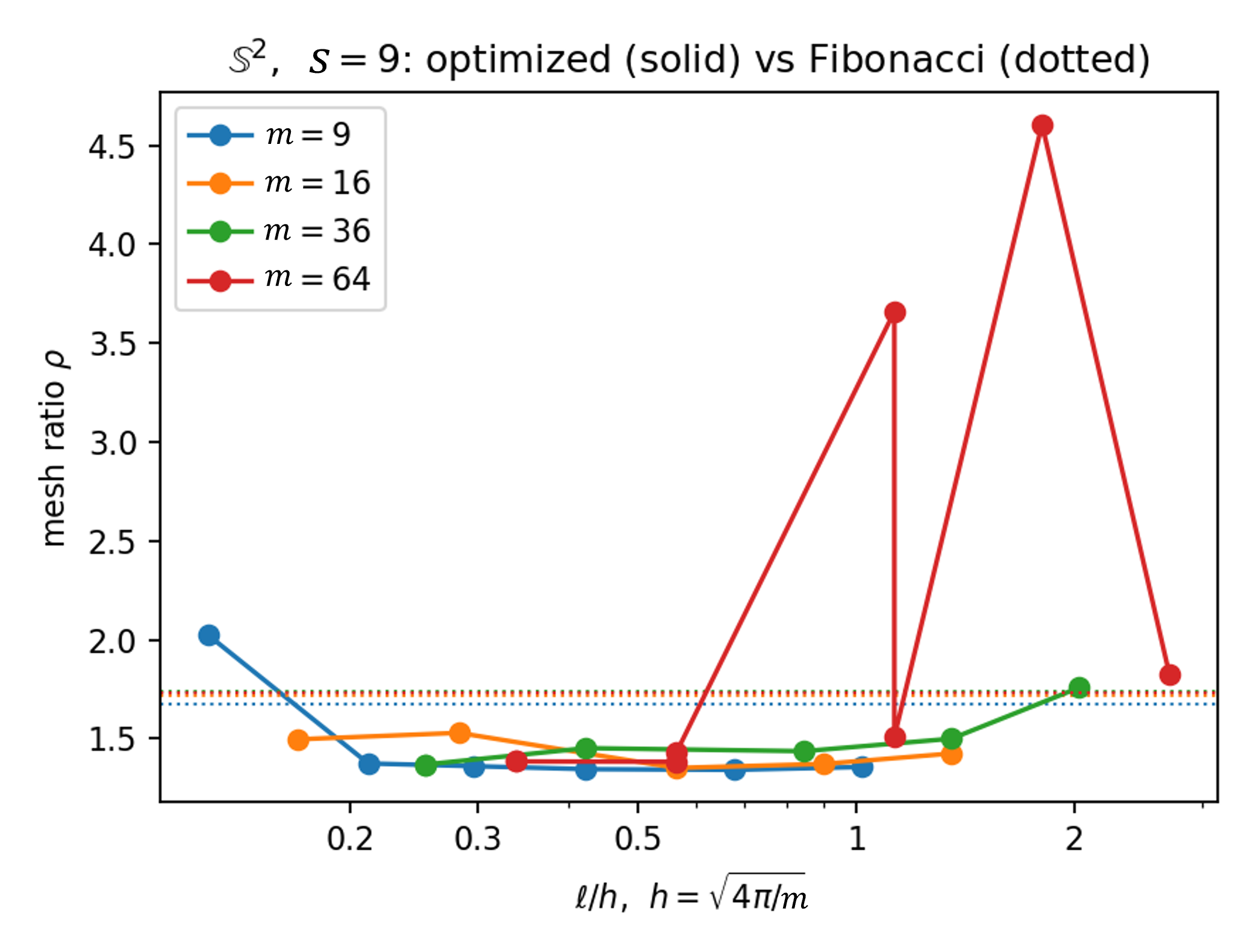}
\caption{Trade-off on $\Sph^2$, $s=9$. Median mesh ratio $\rho$ of optimized landmarks (solid) and the Fibonacci lattice with the same $m$ (dotted, same colour)}
\label{fig:tradeoff_sphere}
\end{figure}

The results from the circle extend to the sphere, where evenness was measured by the \emph{mesh ratio} $\rho=\delta_Z/q_Z$, where $\delta_Z$ is the covering radius and $q_Z$ is half the minimum separation.\footnote{Some references define the mesh ratio as the covering radius over the minimum separation.} We use $\delta_Z$ for the covering radius to avoid confusion with the spacing $h$. In all but a few of $(m,\ell)$ combinations shown in Figure~\ref{fig:tradeoff_sphere}
the optimized landmarks have a lower mesh ratio than the Fibonacci lattice, and in all cases they reach a lower objective (Figure~\ref{fig:tradeoff_sphere}). The advantage in $F$ is largest at moderate-to-large $\ell/h$, where the objective is far from its flat limit. There $F_{\rm opt}/F_{\rm Fib}$ reaches $0.036$ at $(m,\ell)=(64,0.5)$ and $0.015$ at $(36,1.2)$. For $\ell/h\lesssim0.3$ the two objectives agree to two or three digits, consistent with the flat-landscape limit $F\approx s-m(1-\tau_s)$. Even in this flat regime the optimized sets remain well spread ($\rho\le1.53$), because the weak residual repulsion between neighboring translates is still enough to separate landmarks, although the result is then less clearly tied to the eigenspace.

The single exception is $(m,\ell)=(36,1.2)$, where $\ell/h=2.0$, $F_{\rm opt}=2.7\times10^{-7}$, and $\rho_{\rm opt}=1.76$. This is the sphere analogue of the upper failure on the circle, in which the objective is so small that absolute stopping criteria are triggered before the minimum is reached. For the values of $m$ tested, the useful window on the sphere is therefore roughly $0.3\lesssim\ell/h\lesssim1.4$, where the objective is both informative and well resolved. 

\section{Results on The Circle}\label{sec:circle}
To first establish a baseline for accuracy, we use the exact eigenspace of
Proposition~\ref{prop:circle} so that the objective is free of sampling noise.
From $20$ random starts per case, we ask whether the optimizer reliably finds
the equispaced configuration, and whether the attained $F$ agrees with the
closed-form value of Theorem~\ref{thm:alias}. We then turn to a data-driven
experiment, in which the eigenvectors are computed from $N=10{,}000$ uniform
samples on $\Sph^1$ with $\ell=0.8$.

\subsection{Exact-Eigenspace Study}

Table~\ref{tab:circle_multi} reports $20$ random starts per case with the analytic objective at $\ell=0.8$. In every case with $m\le10$ all $20$ runs reached the equispaced configuration, the best $F$ found matched the value from Theorem~\ref{thm:alias}, and no start found a lower value. The median gap deviations of order $10^{-9}$ to $10^{-8}$ show that the equispacing is exact to the optimizer's tolerance. The Hessian at the equispaced configuration is positive definite modulo the rotation direction in all of these cases. The single failure is $(s,m)=(9,17)$, where none of the runs achieves equispacing. In that case, the leakage $7.8\times10^{-10}$ lies below the level of
about $10^{-9}$ at which random starts fail to reach equispacing in the trade-off study
(Section~\ref{sec:conditioning}), and the Gram matrix is poorly conditioned with
$\kappa(K_{ZZ})\approx4.3\times10^5$. Here $h=2\pi/17$, so $\ell/h\approx2.2$, which lies above the window identified in Section~\ref{sec:conditioning} for $m$ of this size. This is the behavior that
Section~\ref{sec:conditioning} predicts, and Remark~\ref{rem:smallF} suggests how it might be avoided.

\begin{table}[h!]
\centering
\small
\resizebox{\textwidth}{!}{\begin{tabular}{@{}rrcccccc@{}}
\toprule
$s$ & $m$ & $F$ equispaced (Thm.~\ref{thm:alias}) & best $F$ found & runs at equispaced & median gap dev. & Hessian $\eta_2$ & $\kappa(K_{ZZ})$\\
\midrule
3 & 3 & 0.643775 & 0.643775 & 20/20 & $6\times10^{-10}$ & 0.588 & 1.3\\
5 & 5 & 0.434518 & 0.434518 & 20/20 & $4\times10^{-10}$ & 1.072 & 3.7\\
7 & 7 & 0.333844 & 0.333844 & 20/20 & $1\times10^{-9}$ & 0.454 & 15.6\\
9 & 9 & 0.273304 & 0.273304 & 20/20 & $1\times10^{-8}$ & 0.092 & 84.7\\
5 & 10 & $2.0428\times10^{-5}$ & $2.0428\times10^{-5}$ & 20/20 & $3\times10^{-9}$ & $1.7\times10^{-3}$ & 319\\
9 & 17 & $7.818\times10^{-10}$ & $7.867\times10^{-10}$ & 0/20 & 0.59 & $2.0\times10^{-7}$ & $4.3\times10^{5}$\\
\bottomrule
\end{tabular}}
\caption{Circle, $\ell=0.8$, exact eigenspace. ``Hessian $\eta_2$ '' is the smallest nonzero eigenvalue of the Hessian of $F$ at the equispaced configuration, whose eigenvalues we denote $\eta_1\le\eta_2\le\cdots$ to distinguish them from the kernel eigenvalues ($\eta_1=0$ is the rotation). No start found a value below the equispaced one. For $(s,m)=(9,17)$ the objective is below the level of about $10^{-9}$ at
which the trade-off study begins to fail, the Gram matrix is badly conditioned, and the runs stall away from equispacing.}
\label{tab:circle_multi}
\end{table}

Table~\ref{tab:tight} complements this study by quantifying Theorem~\ref{thm:tight}. The equispaced sets on the circle attain the maximum possible $\sigma_{\min}(\Phi)^2=s$ and are therefore exactly tight, and the same holds for the tetrahedron on the sphere. For the larger sphere cases the optimized sets are far closer to tight than Fibonacci lattices, even though exact tightness is impossible there.

\begin{table}[h!]
\centering
\resizebox{\textwidth}{!}{\begin{tabular}{@{}lrcccc@{}}
\toprule
& $s$ & $\sigma_{\min}(\Phi)^2$ & max possible & bound $s\tau_s/(\lambda_s\sigma_{\min}^2)$ & $F$\\
\midrule
$\Sph^1$ equispaced & 3 & 3.00 & 3 & 0.928 & 0.644\\
$\Sph^1$ equispaced & 5 & 5.00 & 5 & 0.607 & 0.435\\
$\Sph^1$ equispaced & 7 & 7.00 & 7 & 0.446 & 0.334\\
$\Sph^1$ equispaced & 9 & 9.00 & 9 & 0.350 & 0.273\\
$\Sph^2$ optimized (tetrahedron) & 4 & 4.00 & 4 & 2.020 & 1.214\\
$\Sph^2$ optimized & 9 & 5.70 (Fib.\ 3.03) & 9 & 2.998 & 1.360\\
$\Sph^2$ optimized & 16 & 10.21 (Fib.\ 1.21) & 16 & 2.865 & 1.451\\
\bottomrule
\end{tabular}} \caption{Theorem~\ref{thm:tight} in numbers ($\ell=0.8$). The optimized sets are exactly tight where tightness is possible (equispaced on $\Sph^1$, tetrahedron on $\Sph^2$) and far closer to tight than Fibonacci lattices otherwise.}
\label{tab:tight}
\end{table}

\subsection{Data-Driven Experiment} 




The experiment uses $N=10{,}000$ uniform points on $\Sph^1$ with $\ell=0.8$. For $s=m=5$ with random initializations (and one greedy initialization) all runs reach the same minimum $F=0.4215$. This suggests a single basin, up to the rotation and permutation symmetries of the circle. The angular spacings of the optimized landmarks sum to $2\pi\approx6.2832$ and sit within about $1\%$ of $2\pi/5=1.2566$, with a largest deviation of $0.0124$ rad. These deviations are small but nonzero because the empirical measure $\mu_N$ is not rotation invariant. The data-driven objective therefore lacks the symmetry that makes the equispaced configuration a critical point, and its minimizer is a slight perturbation of an equispaced set. (The same lack of invariance means that rotating all landmarks together changes $F$ slightly, which fixes the overall rotation of the optimum) With the exact eigenspace the deviations vanish to $10^{-9}$ (Table~\ref{tab:circle_multi}).

We also compute the principal angles directly. With $A=LL^\top$, the singular
values of $Y=L^{-1}\Psi$ are the cosines of the principal angles between
$V_5$ and $U_{\hat Z}$, where $\Psi=k(\hat Z,X)H$. For the optimal landmarks
$\hat Z$ from the optimizer the five squared cosines sum to
\[
0.99695+0.98379+0.98341+0.81130+0.80305=4.57850,
\]
so $F=5-4.57850=0.42150$, which equals the optimizer's best $F$ to the printed precision. The Cholesky objective and the SVD diagnostic therefore agree. The leakage $1-\cos^2\theta$ in each direction also matches the Fourier picture. It is approximately $0.003$ for the constant mode, $0.016$ and $0.017$ for the $k=1$ pair, and $0.189$ and $0.197$ for the $k=2$ pair. The leakage is largest for $k=2$ because the landmark translates carry the $k=\pm3$ content (the $\lambda_6,\lambda_7$ pair) aliased onto $k=\mp2$, and those $k=3$ modes are the closest ones in the tail. This is exactly the mechanism of Theorem~\ref{thm:alias}, which gives a per-mode leakage of $\sin^2\theta_k=\tan^2\theta_k/(1+\tan^2\theta_k)$ with $\tan^2\theta_k=\sum_{j\ne0}\beta_{|k+5j|}/\beta_{|k|}$. With the exact eigenvalues this yields $0.0031$, $0.0168$, and $0.1989$ for $k=0,1,2$, in line with the measured values. By mass conservation \eqref{eq:massconservation}, the same $0.4215$ is also the total amount of the tail $i>5$ that has leaked into $U_Z$.

The sandwich bound can be checked directly on these numbers. Since $k(x,x)=1$ we have $\sum_i\lambda_i=1$, and therefore
\[
\sum_{i>5}\lambda_i=1-(0.35847+0.21934+0.21888+0.07916+0.07683)=1-0.95268=0.04732.
\]
The mean of $P_{\hat Z}^2$ over the data is $0.07810$, which estimates $\int P_{\hat Z}^2\dd$, so the excess is $E=0.07810-0.04732=0.03078$. The sandwich bound $(\lambda_5-\lambda_6)\,F\le E\le\lambda_1F$ gives a lower bound of $(0.07683-0.01955)(0.4215)=0.0241$ and an upper bound of $0.3585\times0.4215=0.1511$. The observed value $E=0.0308$ lies inside $[0.0241,\,0.1511]$ and near its lower end. This is expected, because most of the leakage sits on the small eigenvalues (the $k=2$ modes), so the weighted excess is small compared to the worst case where all the mass lies on $\lambda_1$.

\begin{figure}[h!]
\centering
\includegraphics[width=\linewidth]{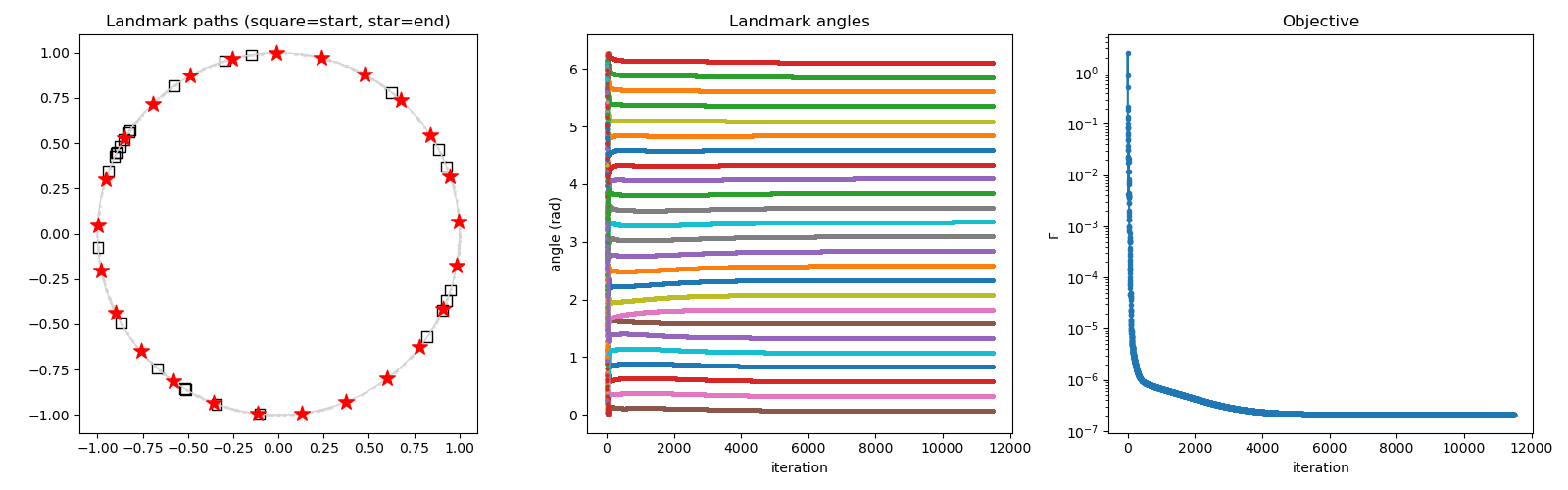}
\caption{Circle, $s=9$, $m=25$, $\ell=0.30$ ($\ell/h\approx1.19$), data-driven eigenvectors from $N=10{,}000$ uniform samples.}
\label{fig:circle_datadriven}
\end{figure}

The experiments with more landmarks than eigenvectors probe the regime of Theorem~\ref{thm:alias}. With $s=5$ and $m=10$, the leakage falls to $F=2.0613\times10^{-5}$ very quickly and equispacing is maintained, in close agreement with the value $2.0428\times10^{-5}$ predicted for the exact eigenspace. Figure~\ref{fig:circle_datadriven} shows a larger data-driven example, with $s=9$, $m=25$ and $\ell=0.30$, so that $h=2\pi/25$ and $\ell/h\approx1.19$, inside the window of Section~\ref{sec:conditioning}.

\section{Results on the Sphere}\label{sec:sphere}
The sphere experiments test whether the picture established on the circle survives in a setting where exact tightness is mostly impossible~\cite{BannaiDamerell1979}. We examine which configurations the optimizer finds, whether the symmetric critical points predicted by Theorem~\ref{thm:symmetry} are minima, and how the sandwich bound behaves. Unless stated otherwise, the objective uses the exact truncated kernel $K^{(s)}$ of Proposition~\ref{prop:sphere}.

Landmarks are parametrized by unconstrained vectors $v_j\in\R^3$ with $z_j=v_j/\|v_j\|$, and they are renormalized after every optimizer step. Since $F$ depends on $v_j$
only through its direction, the gradient with respect to $v_j$ is tangent to the sphere at $z_j$, and renormalization changes neither $F$ nor its gradient direction. We used Adam with learning rate $0.03$, followed by L-BFGS with a strong Wolfe line search, a gradient-and-patience stopping rule, jitter $\varepsilon=10^{-13}$, and $\ell=0.8$.

Table~\ref{tab:poly} evaluates $F$ at the Platonic solids. As Theorem~\ref{thm:symmetry} predicts, all are critical points to machine precision. The signature of the Euclidean Hessian of $Z\mapsto F(Z/\|Z\|)$ counts negative, zero and positive eigenvalues. The zero eigenvalues are always $m+3$ in number, corresponding to $m$ radial directions plus $3$ rotations.

\begin{table}[h!]
\centering
\begin{tabular}{@{}lrrccl@{}}
\toprule
configuration & $m$ & $s$ & $F$ & $\|\nabla F\|$ & Hessian $(-,0,+)$\\
\midrule
tetrahedron & 4 & 4 & 1.213989 & $5\times10^{-17}$ & $(0,7,5)$ \textbf{minimum}\\
octahedron & 6 & 4 & 0.421817 & $1\times10^{-16}$ & $(0,9,9)$ \textbf{minimum}\\
cube & 8 & 4 & 0.209985 & $4\times10^{-16}$ & $(2,11,11)$ saddle\\
icosahedron & 12 & 4 & 0.012389 & $9\times10^{-16}$ & $(12,15,9)$ saddle\\
icosahedron & 12 & 9 & 0.324244 & $1\times10^{-15}$ & $(0,15,21)$ \textbf{minimum}\\
dodecahedron & 20 & 9 & 0.104477 & $3\times10^{-14}$ & $(28,23,9)$ saddle\\
dodecahedron & 20 & 16 & 0.728297 & $2\times10^{-14}$ & $(3,23,34)$ saddle\\
\bottomrule
\end{tabular}
\caption{Platonic solids as critical points of $F$ on $\Sph^2$ ($\ell=0.8$).}
\label{tab:poly}
\end{table}

\begin{figure}[h!]
\centering
\includegraphics[width=\linewidth]{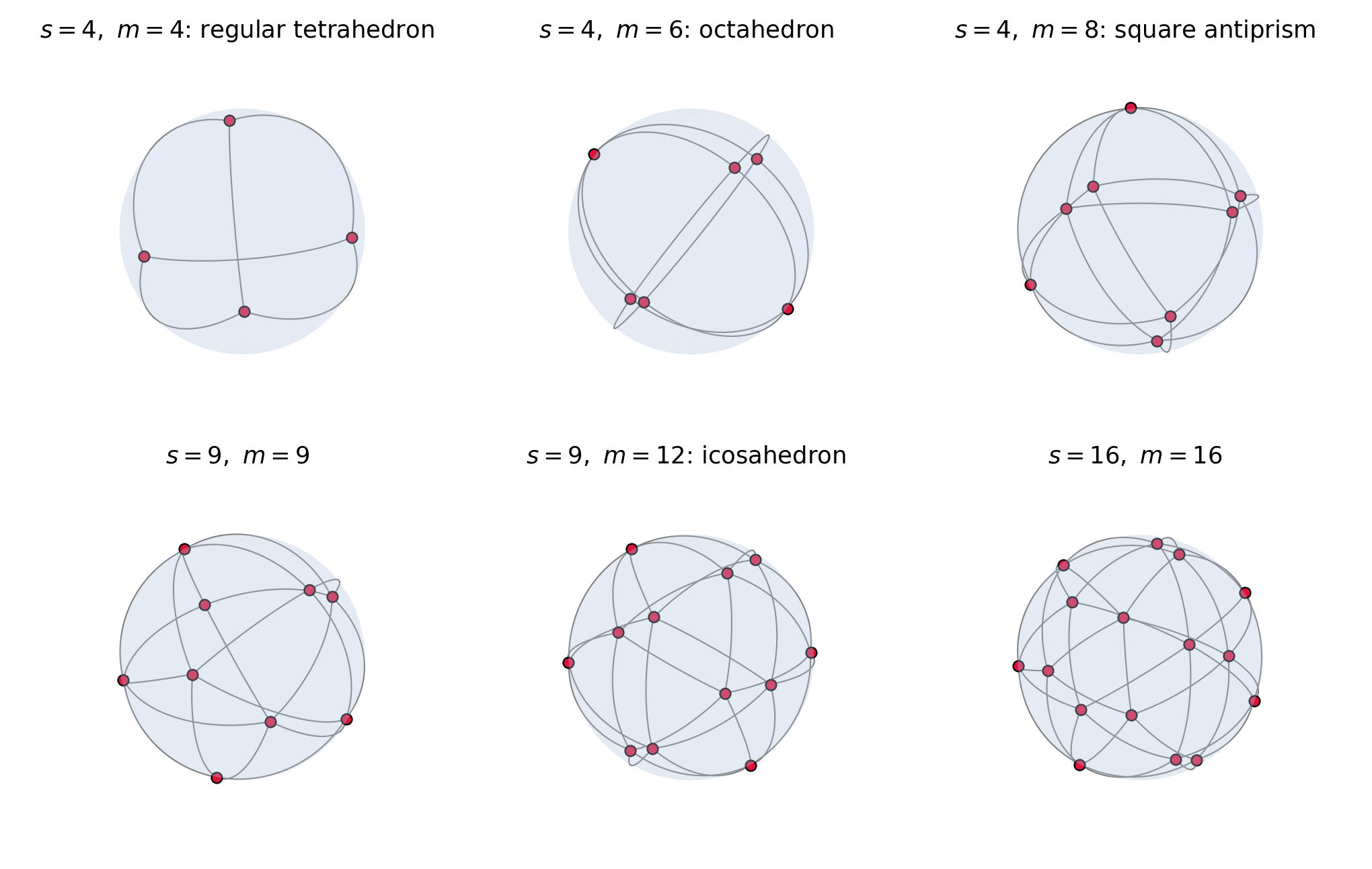}
\caption{Best landmark configurations on $\Sph^2$ from ten random starts ($\ell=0.8$, exact eigenspace), drawn with the geodesic edges of their convex hull.}
\label{fig:sphere_configs}
\end{figure}

Table~\ref{tab:sphere_multi} summarizes ten random starts per case. Figure~\ref{fig:sphere_configs} shows the best configurations. For $m=4$, $6$ and $12$ the identification is rigorous up to numerical tolerance, through the Tammes problem of maximizing the minimum separation of $m$ points \cite{Tammes1930}. Fejes T\'oth's bound \cite{FejesToth1943} is attained for $m=4,6,12$, and the maximizers are unique up to isometry. We therefore claim to recover the regular tetrahedron, octahedron and icosahedron, with separations $\arccos(-\tfrac13)=109.4712^\circ$, $90^\circ$ and $\arctan 2=63.4349^\circ$ (for uniqueness at $m=12$ see also \cite{Danzer1986}). In the best run (which all ten starts reach) the nearest-neighbor angle of \emph{every} landmark equals these values to within $10^{-7}$ degrees. An independent test also shows that the best $F$ equals $F$ at the explicitly constructed polyhedron (Table~\ref{tab:poly}) to about $15$ significant digits.

\begin{table}[h!]
\centering
\small
\resizebox{\textwidth}{!}{\begin{tabular}{@{}rrccccccl@{}}
\toprule
$s$ & $m$ & best $F$ & runs at best & $F$ Fibonacci & min sep.\ & mesh ratio & mesh Fib.\ & identified\\
\midrule
4 & 4 & 1.213989 & 10/10 & 1.2553 & $109.47^\circ$ & 1.29 & 1.57 & regular tetrahedron\\
4 & 6 & 0.421817 & 10/10 & 0.5488 & $90.00^\circ$ & 1.21 & 1.60 & octahedron\\
4 & 8 & 0.169303 & 10/10 & 0.2114 & $71.90^\circ$ & 1.56 & 1.65 & square antiprism\\
9 & 9 & 1.359921 & 10/10 & 1.5085 & $69.37^\circ$ & 1.34 & 1.67 & tricapped trigonal prism ($D_{3h}$)\\
9 & 12 & 0.324244 & 10/10 & 0.5756 & $63.43^\circ$ & 1.18 & 1.70 & icosahedron\\
16 & 16 & 1.451253 & 10/10 & 1.8719 & $49.39^\circ$ & 1.34 & 1.72 & isotropic ($2$-design)\\
16 & 31 & $7.2448\times10^{-3}$ & 9/10 & $2.265\times10^{-2}$ & $35.32^\circ$ & 1.37 & 1.74 & --\\
\bottomrule
\end{tabular}}
\caption{10 Random Starts on the Sphere, $\ell=0.8$, exact eigenspace. ``Runs at best'' counts starts within relative $10^{-6}$ of the best value achieved. Every optimized set beats the Fibonacci lattice in both $F$ and mesh ratio. The tetrahedron, octahedron and icosahedron are recovered exactly: their $F$ values match Table~\ref{tab:poly} to all printed digits, and all nearest-neighbor angles agree to $10^{-7}$ degrees.}
\label{tab:sphere_multi}
\end{table}

To further verify the tightness claims, we certify the design strength against the number of points. A finite set $Z=\{z_1,\dots,z_m\}\subset\Sph^2$ is a spherical $t$-design if
the equal-weight rule integrates every polynomial of degree at most $t$
exactly,
\[
\frac{1}{m}\sum_{j=1}^m p(z_j)=\int_{\Sph^2}p\,d\mu\qquad\text{for all }p\in\Pi_t ,
\]
and its \emph{strength} is the largest such $t$ \cite{DelsarteGoethalsSeidel1977}.
Since the restriction of $\Pi_t$ to $\Sph^2$ is $\bigoplus_{n\le t}\mathcal H_n$,
where $\mathcal H_n$ is the space of spherical harmonics of degree $n$, and
since every $Y\in\mathcal H_n$ with $n\ge1$ has zero mean, the definition
reduces to a condition on each degree separately. For $n\ge1$ define the
\emph{degree-$n$ moment}
\begin{equation}\label{eq:moment}
M_n(Z):=\Big(\frac{1}{2n+1}\sum_{k=1}^{2n+1}\Big|\frac{1}{m}\sum_{j=1}^m\varphi_{nk}(z_j)\Big|^2\Big)^{1/2}
=\Big(\frac{1}{m^2}\sum_{j,l=1}^mP_n(z_j\cdot z_l)\Big)^{1/2},
\end{equation}
where $\{\varphi_{nk}\}_k$ is an $L^2(\mu)$-orthonormal basis of
$\mathcal H_n$. The second expression follows from the addition theorem used
in Proposition~\ref{prop:sphere}. It shows that $M_n$ does not depend on the
basis and is invariant under rotations, reflections and relabelling of $Z$.
Then
\[
Z\text{ is a $t$-design}\iff M_n(Z)=0\ \text{ for } 1\le n\le t ,
\]
so the strength is the largest $t$ for which $M_1,\dots,M_t$ all vanish.
Two special cases are worth noting. First, $M_1(Z)=\|\frac{1}{m}\sum_jz_j\|$ is
the length of the centroid. Second, $M_1=M_2=0$ is equivalent to a zero
centroid together with the isotropy condition $\sum_jz_jz_j^\top=\frac{m}{3}I$,
which is the $2$-design criterion used for $(16,16)$ below.

We evaluated \eqref{eq:moment} through the harmonic averages (first form),
not the Legendre double sum. The double sum equals $M_n^2$, and its roundoff
of about $10^{-16}$ would hide any moment below about $10^{-8}$. The harmonic
averages are linear in the data and resolve moments down to about $10^{-12}$.
Table~\ref{tab:moments} reports $M_n$ for $1\le n\le7$ at the best run of
each case. Every entry is either
at most $1.7\times10^{-9}$, which is the size expected from the stopping rule
(gradient below $3\times10^{-9}$), or at least $2.3\times10^{-3}$.
Nothing lies in between, a gap of six orders of magnitude. We call $M_n$
vanishing when $M_n<10^{-8}$. Any threshold in the gap gives the same
strengths:
\[
t(4,4)=2,\quad t(4,6)=3,\quad t(4,8)=1,\quad t(9,9)=1,\quad t(9,12)=5,\quad t(16,16)=2,
\]
where $t(s,m)$ denotes the strength of the best configuration for the pair $(s,m)$. The $(16,31)$ minimizer is not a $1$-design. Its centroid has length $M_1=1.0\times10^{-2}$, and all of $M_1,\dots,M_7$ lie between $2.3\times10^{-3}$ and $3.1\times10^{-2}$, so it is not a design of any strength. The strengths
of the tetrahedron, octahedron and icosahedron agree with their known values
\cite{DelsarteGoethalsSeidel1977}.

\begin{table}[h]
\centering
\small
\resizebox{\textwidth}{!}{\begin{tabular}{@{}rr*{7}{c}c@{}}
\toprule
$s$ & $m$ & $M_1$ & $M_2$ & $M_3$ & $M_4$ & $M_5$ & $M_6$ & $M_7$ & strength $t$\\
\midrule
4 & 4 & $4.1{\times}10^{-10}$ & $5.8{\times}10^{-10}$ & $0.75$ & $0.51$ & $1.4{\times}10^{-9}$ & $0.63$ & $0.61$ & 2\\
4 & 6 & $1.7{\times}10^{-12}$ & $3.2{\times}10^{-11}$ & $6.6{\times}10^{-11}$ & $0.76$ & $8.8{\times}10^{-11}$ & $0.35$ & $1.3{\times}10^{-10}$ & 3\\
4 & 8 & $3.5{\times}10^{-10}$ & $3.4{\times}10^{-2}$ & $3.0{\times}10^{-10}$ & $0.37$ & $0.59$ & $0.26$ & $0.26$ & 1\\
9 & 9 & $3.5{\times}10^{-11}$ & $3.3{\times}10^{-3}$ & $7.9{\times}10^{-2}$ & $0.14$ & $0.61$ & $0.35$ & $5.8{\times}10^{-2}$ & 1\\
9 & 12 & $1.3{\times}10^{-10}$ & $4.2{\times}10^{-10}$ & $6.2{\times}10^{-10}$ & $8.4{\times}10^{-10}$ & $7.0{\times}10^{-10}$ & $0.66$ & $1.0{\times}10^{-9}$ & 5\\
16 & 16 & $1.9{\times}10^{-11}$ & $1.8{\times}10^{-10}$ & $6.3{\times}10^{-2}$ & $2.2{\times}10^{-2}$ & $1.7{\times}10^{-9}$ & $0.19$ & $0.52$ & 2\\
16 & 31 & $1.0{\times}10^{-2}$ & $1.1{\times}10^{-2}$ & $6.3{\times}10^{-3}$ & $2.3{\times}10^{-3}$ & $1.6{\times}10^{-2}$ & $3.1{\times}10^{-2}$ & $2.6{\times}10^{-2}$ & --\\
\bottomrule
\end{tabular}}
\caption{Degree-$n$ moments \eqref{eq:moment} of the best configuration in
each case ($\ell=0.8$, exact eigenspace). Every entry is either below
$1.7\times10^{-9}$ or above $2.3\times10^{-3}$. The strength is the largest
$t$ with $M_1,\dots,M_t<10^{-8}$.}
\label{tab:moments}
\end{table}

The geometry of the minimizers that are not Platonic solids is also informative. For $s=4$ and $m=8$ the optimizer finds a square antiprism rather than the cube. Its $16$ edges come in two lengths, $71.90^\circ$ and $79.82^\circ$, so it is not the equilateral Tammes antiprism~\cite{SchutteVanDerWaerden1951}, but it has the same combinatorial type. This identification rests on the edge pattern only, we did not test for the $D_{4d}$ symmetry of the antiprism. The optimizer for $s=m=9$ has twelve nearest-neighbor edges of equal length $69.37^\circ$ and forms a tricapped trigonal prism, consisting of two equilateral triangles at heights $\pm0.7095$ and an equatorial triangle rotated by $60^\circ$.  This is the combinatorial type of the $9$-point Tammes solution \cite{SchutteVanDerWaerden1951}, but with a smaller separation ($69.37^\circ$ against $\arccos(1/3)=70.53^\circ$), so it is not the Tammes configuration. We checked that it is invariant under a $3$-fold rotation and the horizontal mirror to within $5\times10^{-9}$, so its symmetry group contains $D_{3h}$, and its inertia tensor is isotropic to within $1\%$. The $m=16,s=16$ optimizer has zero centroid and an exactly isotropic inertia tensor $\sum_jz_jz_j^\top=\frac{16}{3}I$. These two conditions are equivalent to being a spherical $2$-design \cite{DelsarteGoethalsSeidel1977}. The configuration has six pairs at the tetrahedral angle $109.47^\circ$, which is suggestive of tetrahedral symmetry, but we did not test for it. The $(16,31)$ optimizer has a centroid of length $1.0\times10^{-2}$, and its moments $M_1,\dots,M_7$ are all between $2.3\times10^{-3}$ and $3.1\times10^{-2}$ (Table~\ref{tab:moments}). It is therefore well spread but is not a spherical design, even approximately at the level of the other cases.

\begin{figure}[h!]
\centering
\includegraphics[width=0.9\textwidth]{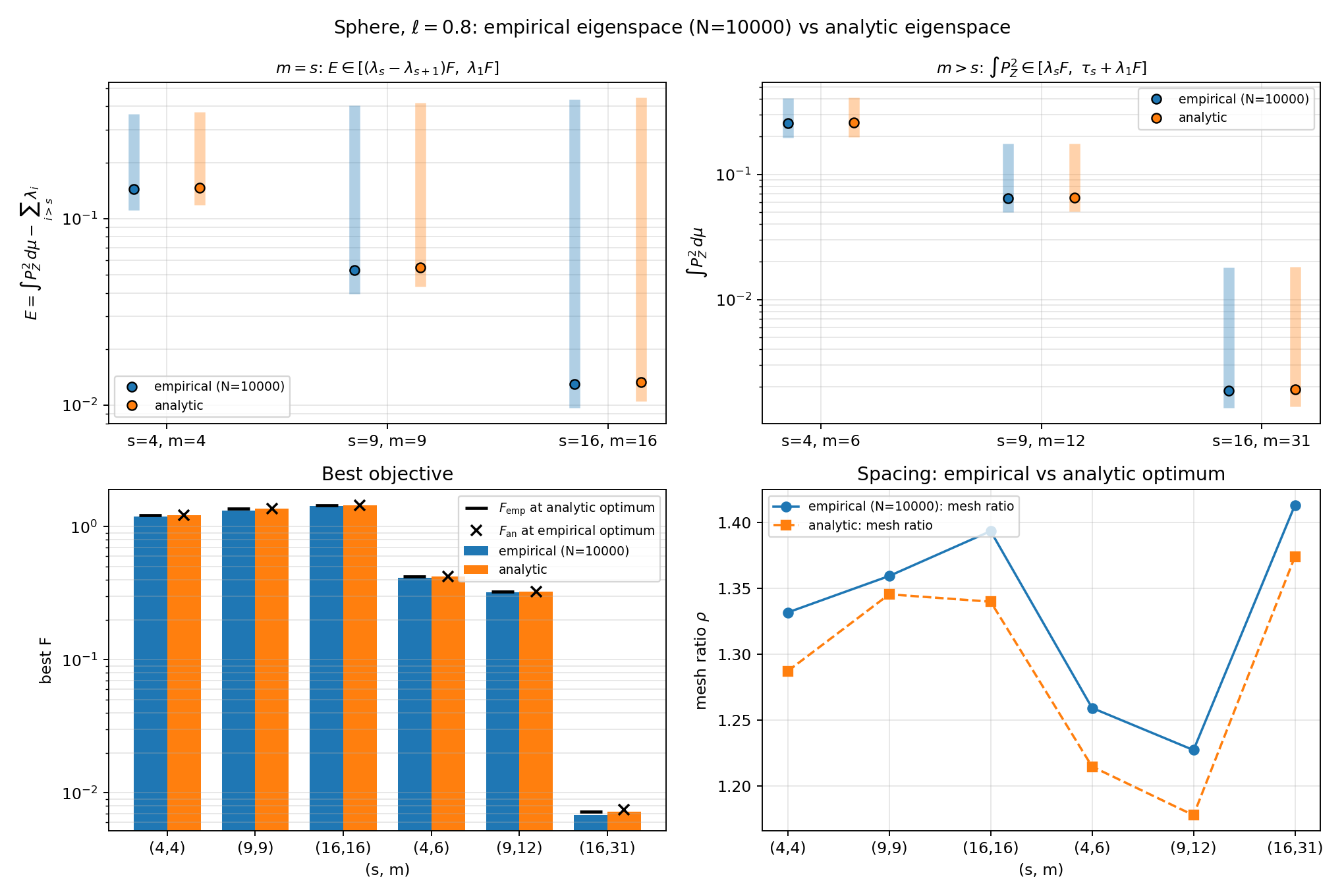}
\includegraphics[width=0.8\linewidth]{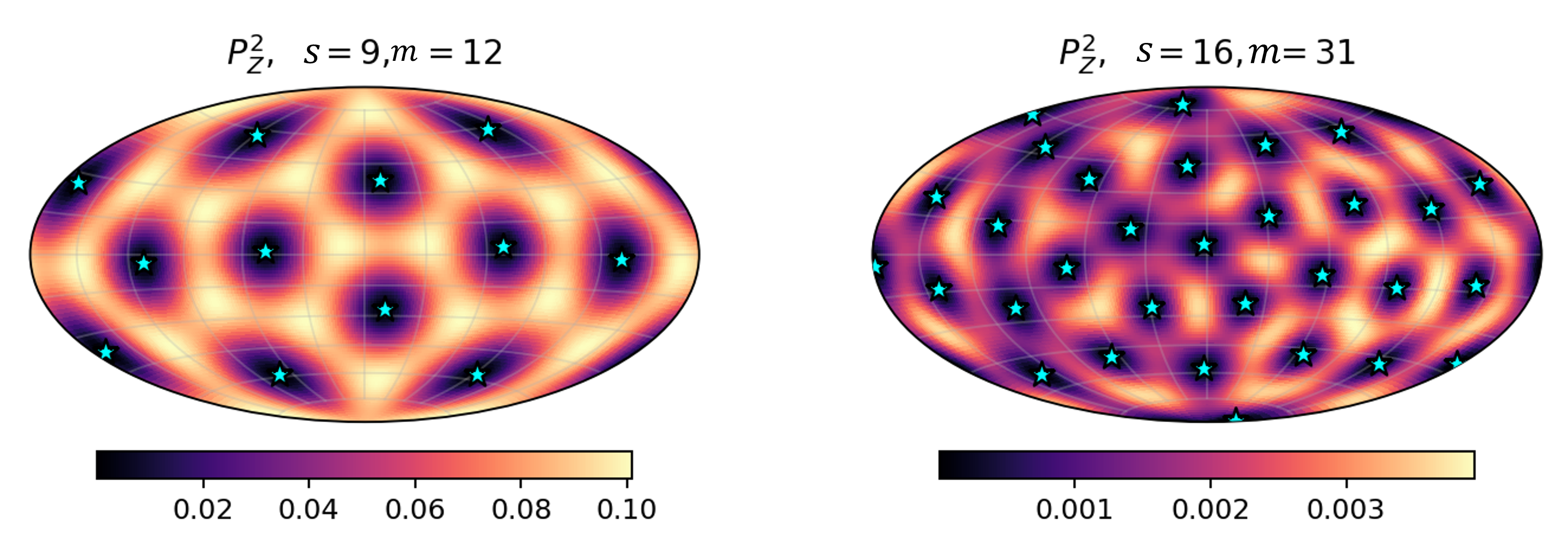}
\caption{Sandwich bounds and optimal landmarks on $\Sph^2$ ($\ell=0.8$) with the empirical eigenspace from $N=10{,}000$ uniform samples (blue) and the exact eigenspace (orange). Each case uses the best of ten random starts.
Top left: for $s=m$, the excess $E=\int P_Z^2\,d\mu-\sum_{i>s}\lambda_i$ (dots) and the interval $[(\lambda_s-\lambda_{s+1})F,\ \lambda_1F]$ of Theorem~\ref{thm:sandwich} (bars).
Top right: for $m>s$, $\int P_Z^2\,d\mu$ and the interval $[\lambda_sF,\ \tau_s+\lambda_1F]$ of Corollary~\ref{cor:sgeqm}.
Bottom left: best objective value under each eigenspace (bars), with cross-evaluations. Dashes give the empirical objective at the analytic optimum and crosses the analytic objective at the empirical optimum.
Bottom right: mesh ratio of the empirical and analytic optima (left axis), and the largest difference between their sorted pairwise geodesic angles (grey bars, right axis, degrees). The angle difference does not depend on rotation, reflection or point labelling.}
\label{fig:sandwich_emp}
\end{figure}

To test whether the conclusions above depend on having the exact eigenspace, we repeated the sphere experiments with the data-driven eigenspace used in practice. We drew $N=10{,}000$ uniform points $X\subset\Sph^2$, computed the full eigendecomposition of the Gaussian kernel matrix ($\ell=0.8$), and took the empirical Mercer eigenvalues and the truncated kernel. With the completed basis of Section~\ref{sec:setting}, which accounts for the components of off-data translates in $\Hcal_0$, Theorem~\ref{thm:sandwich} and Corollary~\ref{cor:sgeqm} hold exactly for $\mu_N$, and we can compare the sandwich bounds. We used the same ten random starts and the same optimizer as for Table~\ref{tab:sphere_multi}, kept the best run, and compared it with the best analytic configuration from that table. Table~\ref{tab:sandwich_emp} and Figure~\ref{fig:sandwich_emp} summarize the results.

For every $(s,m)$ pair, the integrated power function at the best empirical configuration lies at the lower end of the sandwich bound,
 and the best empirical $F$ values are $1$--$6\%$ below the analytic ones. The cross-evaluations show that the analytic objective is barely affected by the perturbation because the analytic optima are nondegenerate minima, so a small displacement changes $F_{\rm an}$ only at second order. The empirical objective lacks the rotational symmetry of the exact problem, and so it pulls the optimum away from the symmetric configuration at first order.


The empirical optima are well spread in every case, with mesh ratios between $1.23$ and $1.41$, well below the Fibonacci values in Table~\ref{tab:sphere_multi}. They are, however, slightly worse than the exact optima, with mesh ratios higher by $0.014$--$0.054$ and minimum separations smaller by $0.6^\circ$--$2.8^\circ$. The sorted pairwise angles of the two optima differ by at most $1.6^\circ$--$4.4^\circ$, with root-mean-square differences of $0.6^\circ$--$1.4^\circ$. The empirical optima are therefore recognizable distortions of the symmetric configurations rather than different configurations. For instance, the nearest-neighbor angles of the empirical $(4,4)$ optimum range over $106.6^\circ$--$109.9^\circ$ instead of the tetrahedral $109.47^\circ$, and those of the $(9,12)$ optimum range over $61.7^\circ$--$63.2^\circ$ instead of the icosahedral $63.43^\circ$.

Sampling breaks the rotational symmetry and splits the degenerate eigenvalues, which perturbs the optimal landmarks by a few degrees. We expect this perturbation to shrink as $N$ grows, since the eigenvalue splitting does, but we did not vary $N$ here. However the qualitative picture that minimizers of $F$ are evenly spread with power function near the lower end of the sandwich bound, as well as approximating the symmetric configurations predicted by the exact theory, remains unchanged.

\begin{table}[h!]
\centering
\small
\resizebox{\textwidth}{!}{%
\begin{tabular}{@{}rrcccccccc@{}}
\toprule
 & & \multicolumn{2}{c}{best $F$} & \multicolumn{2}{c}{cross-evaluation} & \multicolumn{2}{c}{min.\ separation} & \multicolumn{2}{c}{mesh ratio}\\
\cmidrule(lr){3-4}\cmidrule(lr){5-6}\cmidrule(lr){7-8}\cmidrule(lr){9-10}
$s$ & $m$ & empirical & analytic & $F_{\rm an}(Z_{\rm emp})$ & $F_{\rm emp}(Z_{\rm an})$ & empirical & analytic & empirical & analytic\\
\midrule
4  & 4  & 1.18638 & 1.21399 & 1.21452 & 1.22465 & $106.64^\circ$ & $109.47^\circ$ & 1.33 & 1.29\\
9  & 9  & 1.32081 & 1.35992 & 1.36047 & 1.37020 & $68.47^\circ$  & $69.37^\circ$  & 1.36 & 1.35\\
16 & 16 & 1.42296 & 1.45125 & 1.45235 & 1.43972 & $48.80^\circ$  & $49.39^\circ$  & 1.39 & 1.34\\
4  & 6  & 0.41402 & 0.42182 & 0.42256 & 0.42425 & $87.47^\circ$  & $90.00^\circ$  & 1.26 & 1.21\\
9  & 12 & 0.32039 & 0.32424 & 0.32545 & 0.32600 & $61.68^\circ$  & $63.43^\circ$  & 1.23 & 1.18\\
16 & 31 & $6.796\times10^{-3}$ & $7.245\times10^{-3}$ & $7.443\times10^{-3}$ & $7.223\times10^{-3}$ & $34.61^\circ$ & $35.32^\circ$ & 1.41 & 1.37\\
\bottomrule
\end{tabular}}
\caption{Empirical ($N=10{,}000$) versus exact eigenspace on $\Sph^2$, $\ell=0.8$, best of ten starts. The cross-evaluation columns give the analytic objective at the empirical optimum and the empirical objective at the analytic optimum. Each exceeds the best value of its own objective, so each optimum is specific to its objective.}
\label{tab:sandwich_emp}
\end{table}

\FloatBarrier

\bibliographystyle{unsrt}
\bibliography{references}

\appendix

\section{The circle and the sphere}\label{app:circlesphere}

\subsection{Exact eigenspaces}\label{sec:homog}

The objective $F$ depends on the kernel only through $K_{ZZ}$ and $K^{(s)}_{ZZ}$, so it can be evaluated in closed form whenever the Mercer expansion is known explicitly. This happens on the circle and the sphere. With $\mu$ the normalized uniform measure, the Gaussian kernel $K(x,y)=\exp(-\|x-y\|^2/2\ell^2)$ is invariant under the rotation group, so its eigenspaces are the irreducible representations of that group and its eigenvalues can be written down. Throughout we set $\alpha:=1/\ell^2$, and we index the distinct eigenvalues by frequency (on $\Sph^1$) or degree (on $\Sph^2$), writing them $\beta_k$ or $\beta_n$ to distinguish them from the eigenvalues $\lambda_i$ listed with multiplicity. On the unit sphere $\|x-y\|^2=2-2\,x\cdot y$, and the kernel therefore takes the form $K(x,y)=e^{-\alpha}e^{\alpha\,x\cdot y}$, which is the starting point for both expansions below. We begin with the circle, where the eigenvalues are modified Bessel functions of integer order.

\begin{proposition}[Circle]\label{prop:circle}
On $\Sph^1$ with $\mu=d\theta/2\pi$,\[K(\theta,\phi)=\sum_{k\in\mathbb Z}\beta_{|k|}\,e^{ik(\theta-\phi)},\qquad \beta_k=e^{-\alpha}I_k(\alpha),
\] 
where $I_k$ is the modified Bessel function of the first kind. The $L^2(\mu)$-normalized real eigenfunctions are $1,\sqrt2\cos k\theta,\sqrt2\sin k\theta$, and $\beta_k$ has multiplicity $2$ for $k\ge1$, so that $\lambda_1=\beta_0$ and $\lambda_{2k}=\lambda_{2k+1}=\beta_k$. Gaps occur only at $s=2M+1$. For such $s$,
\[
K^{(s)}(\theta,\phi)=\beta_0+2\sum_{k=1}^M\beta_k\cos k(\theta-\phi),\qquad \tau_s:=\sum_{i>s}\lambda_i=1-\beta_0-2\sum_{k=1}^M\beta_k .
\] 
\end{proposition}
\begin{proof}
This follows from the generating function $e^{\alpha\cos t}=\sum_kI_k(\alpha)e^{ikt}$ with $t=\theta-\phi$, since $x\cdot y=\cos(\theta-\phi)$. The tail formula uses $K(x,x)=1$, so $\sum_i\lambda_i=1$.
\end{proof}

On the sphere the same construction goes through with Legendre polynomials in place of complex exponentials, and the Bessel functions acquire half-integer order. The eigenspaces now have growing dimension, which is why gaps appear at the perfect squares.

\begin{proposition}[Sphere]\label{prop:sphere}
On $\Sph^2$ with $\mu=d\sigma/4\pi$,
\[K(x,y)=\sum_{n\ge0}\beta_n\,(2n+1)\,P_n(x\cdot y),\qquad \beta_n=e^{-\alpha}i_n(\alpha)=e^{-\alpha}\sqrt{\tfrac{\pi}{2\alpha}}\,I_{n+1/2}(\alpha),
\]
where $P_n$ is the Legendre polynomial and $i_n$ the modified spherical Bessel function. The eigenspace of $\beta_n$ is the space of spherical harmonics of degree $n$, of dimension $2n+1$. Gaps occur only at $s=(M+1)^2$, and then \[K^{(s)}(x,y)=\sum_{n=0}^{M}\beta_n(2n+1)P_n(x\cdot y),\qquad \tau_s=1-\sum_{n\le M}(2n+1)\beta_n.
\] 
\end{proposition}
\begin{proof}
Use the Funk--Hecke expansion $e^{\alpha t}=\sum_n(2n+1)i_n(\alpha)P_n(t)$ and the addition theorem $\sum_{k}\varphi_{nk}(x)\varphi_{nk}(y)=(2n+1)P_n(x\cdot y)$ for $L^2(\mu)$-orthonormal harmonics. This is consistent with \cite{MinhNiyogiYao2006}.
\end{proof}

The practical consequence is that, by Remark~\ref{rem:basis}, the objective needs only $K^{(s)}_{ZZ}$ and never an explicit eigenbasis. On $\Sph^1$ this matrix is a short cosine sum, and on $\Sph^2$ it is a Legendre sum that can be evaluated by the three-term recurrence. Both forms are exact, rotation invariant, and compatible with automatic differentiation. No data, no $N\times N$ kernel matrix, and no eigensolver are required, so the only inputs are the bandwidth $\ell$ and the truncation level $s$.

The analytic spectra also agree with what a data-driven computation produces. For $N=10{,}000$ uniform points on $\Sph^1$ with $\ell=0.8$, the empirical values $\lambda_i=\nu_i/N$ match $\beta_k=e^{-\alpha}I_k(\alpha)$ to within $2\%$, and we observed  a similar agreement on $\Sph^2$ for the empirical spectrum used in Section~\ref{sec:sphere}. The empirical eigenvalues within a degenerate pair are split slightly by sampling noise, whereas the analytic pairs are exactly degenerate. That the empirical values converge to the Mercer eigenvalues as $N$ grows follows from \cite{RosascoBelkinDeVito2010}, which justifies using the analytic expressions as the reference in all later comparisons.

\subsection{Proof of the Aliasing Theorem on the Circle}

\begin{proof}[Proof of Theorem~\ref{thm:alias}]

The kernel is $K(\theta,\theta')=\sum_{n\in\mathbb Z}\beta_{|n|}e^{in(\theta-\theta')}$, where $\beta_{|n|}$ is the eigenvalue of the frequency-$n$ mode. Since $\beta_{|n|}$ is nonincreasing in $|n|$, the top-$s$ eigenspace is
$V_s=\Span\{\psi_k:|k|\le M\}$ with $\psi_k=\sqrt{\beta_{|k|}}\,e^{ik\theta}$ and $s=2M+1$.
We use the complex characters in place of the real eigenbasis. This is allowed by Remark~\ref{rem:basis}: $F$ depends only on
$K^{(s)}_{ZZ}=\sum_{|k|\le M}\beta_{|k|}\,e^{ik(z_j-z_l)}$, which is the same for the real and the complex basis of $V_s$.
All sums $\sum_{j}\beta_{|r+jm|}$ converge because $\sum_n\beta_{|n|}=K(x,x)<\infty$.
Take $z_j=2\pi j/m$, $j=0,\dots,m-1$, and define the orthonormal DFT vectors
$(f_r)_j=m^{-1/2}e^{irz_j}$, $r\in\mathbb Z_m$.

Sampling aliases frequencies, \emph{i.e.} $e^{inz_j}=e^{i(n\bmod m)z_j}$, so the sample vector of $e^{in\theta}$ at $Z$ is $\sqrt m\,f_{n\bmod m}$. Therefore
\[
K_{ZZ}=m\sum_{r\in\mathbb Z_m}\Big(\sum_{j\in\mathbb Z}\beta_{|r+jm|}\Big)f_rf_r^*.
\]
So $K_{ZZ}$ is circulant, with eigenvectors $f_r$ and eigenvalues $\kappa_r=m\sum_{j}\beta_{|r+jm|}>0$.

Likewise $K^{(s)}_{ZZ}=m\sum_{|k|\le M}\beta_{|k|}f_{k\bmod m}f_{k\bmod m}^*$.
If $k\ne k'$ with $|k|,|k'|\le M$, then $0<|k-k'|\le 2M<m$, so $k\not\equiv k'\pmod m$.
Because $m\ge s$, each residue contains at most one $|k|\le M$, so $K^{(s)}_{ZZ}$ has eigenvalue $m\beta_{|k|}$ at residue $k$ and $0$ at the other $m-s$ residues.

By \eqref{eq:quadform} with $\psi_k(Z)=\sqrt{m\beta_{|k|}}\,f_k$ and $K_{ZZ}^{-1}f_k=\kappa_k^{-1}f_k$,
\[
\langle\psi_k,\Pi_Z\psi_{k'}\rangle_{\Hcal}=\psi_k(Z)^*K_{ZZ}^{-1}\psi_{k'}(Z)=\delta_{kk'}\,a_k,
\qquad
a_k=\frac{\beta_{|k|}}{\sum_{j\in\mathbb Z}\beta_{|k+jm|}},
\]
since $f_k^*f_{k'}=\delta_{kk'}$ for these residues.
Summing the leakages, $F(Z)=\sum_{|k|\le M}(1-a_k)$, which is the stated formula
because $1-a_k=\sum_{j\ne0}\beta_{|k+jm|}\big/\sum_{j\in\mathbb Z}\beta_{|k+jm|}$.
(Equivalently, $F=s-\tr(K_{ZZ}^{-1}K^{(s)}_{ZZ})=s-\sum_{|k|\le M}m\beta_{|k|}/\kappa_k$.)

By Step 3, the compression $Q_s\Pi_ZQ_s|_{V_s}$ is diagonal in the orthonormal basis $\{\psi_k\}$ of $V_s$, with eigenvalues $a_k$.
These are the $\cos^2\theta_k$, hence $\sin^2\theta_k=1-a_k$, as claimed.
The same conclusion follows from symmetry. $Q_s\Pi_ZQ_s$ commutes with the rotation by $2\pi/m$, which preserves $Z$.
The characters $e^{ik\theta}$ with $|k|\le M$ are eigenvectors of that rotation with eigenvalues $e^{\mp2\pi ik/m}$, which are pairwise distinct because $|k-k'|\le2M<m$.
A commuting operator is diagonal on a simple spectrum, so each character is an eigenvector of $Q_s\Pi_ZQ_s|_{V_s}$.

The $m$ eigenvalues of $K_{ZZ}$ are $\kappa_r$, $r\in\mathbb Z_m$. The factor $m$ cancels in the ratio of the largest to smallest:
\[
\kappa(K_{ZZ})=\frac{\max_r\sum_j\beta_{|r+jm|}}{\min_r\sum_j\beta_{|r+jm|}}.
\]

The residues of $-M,\dots,M$ now exhaust $\mathbb Z_m$, so $\Phi$ is square and Proposition~\ref{thm:factor} applies.
The matrix $R_Z$ has eigenvalues $\tan^2\theta_k$ and is diagonal in the same basis. Since $\tan^2\theta_k=\sin^2\theta_k/\cos^2\theta_k$,
\[
\tan^2\theta_k=\frac{\sum_{j\ne0}\beta_{|k+jm|}/\sum_{j}\beta_{|k+jm|}}{\beta_{|k|}/\sum_{j}\beta_{|k+jm|}}
=\frac{\sum_{j\ne0}\beta_{|k+jm|}}{\beta_{|k|}}.
\]
\end{proof}

\end{document}